%Paper: 9203004
%From: bifet@math.sunysb.edu (Emili Bifet)
%Date: Thu, 19 Mar 92 14:58:20 EST

%On the Abel-Jacobi Map for Divisors of Higher Rank on a Curve,
%Emili Bifet, Franco Ghione and Maurizio Letizia,
%34 pages, AmS-TeX 2.0

\input amstex
\documentstyle{amsppt}
\magnification=\magstep1
\hoffset=.5cm
\voffset=.6cm
\TagsOnRight
\font\got=eufm10
\def\EE{\widetilde E}
\def\gotP{\text{\got P}}
\def\ov{\overline}
\def\np{\par\noindent}

\def\div{\text{\rm Div}^{r,n}_{C/k}}
\def\rank{\text{rank }}
\def\divb{\text{\bf Div}^{r,n}_{C/k}}
\def\x{\text{\bf x}}
\def\X{\text{\bf X}}
\def\Y{\text{\bf Y}}
\def\y{\text{\bf y}}

\def\bfo{\Cal O}
\def\bfG{\text{\bf G}}
\def\bfS{\text{\bf S}}
\def\bfN{\text{\bf N}}
\def\d{\,\text{d}}
\def\longto{\longrightarrow}
\def\Z{\Bbb Z}
\def\F{\Cal F}
\def\L{\Cal L}

\def\U{\Cal U}
\def\E{\Cal E}
\def\D{\Cal D}
\def\B{\Cal B}
\def\C{\Cal C}
\def\F{\Cal F}
\def\N{\Cal N}

\def\H{\Cal H}
\def\V{\Cal V}

\def\sskip{\smallskip}
\topmatter
\title
On the Abel-Jacobi Map for Divisors
of Higher Rank on a Curve
\endtitle
\leftheadtext{E. Bifet, F. Ghione and M. Letizia}
\rightheadtext{On the Abel-Jacobi Map}
\author
Emili Bifet, Franco Ghione and Maurizio Letizia
\endauthor
\address
Department of Mathematics and
Institute for Mathematical Sciences,
State University of New York,
Stony Brook, NY 11794-3651
( E. Bifet )
\endaddress
\email bifet\@math.sunysb.edu\endemail
\address
Dipartimento di Matematica,
Universit\`a di Roma--Tor Vergata,
via della Ricerca Scientifica,
00133 Roma, Italia
( F. Ghione and M. Letizia )
\endaddress

\endtopmatter
\document
\head 1. Introduction \endhead

The aim of this paper is to  present an
algebro-geometric approach to the
study of the geometry of the moduli
space of stable bundles on a
smooth projective curve defined over an algebraically closed field $k$,
of arbitrary characteristic.
This establishes a bridge between the arithmetic approach of [15,17,9]
and the gauge group approach of [2]. It
might also help explain some of the mysterious analogies observed by
Atiyah and Bott, in the form of ``Weil conjectures'' for the
infinite dimensional algebraic variety associated to the ind-variety
considered below.

One of the basic ideas is to consider a notion of divisor of
higher rank and a suitable Abel-Jacobi map generalizing the classical
notions in rank one.
Let $\bfo_{C}$ be the structure sheaf of the curve $C$ and
let $K$  be its field of
rational functions, considered as a constant $\bfo_{C}$-module.

We define a divisor of rank $r$ and degree $n$, an $(r,n)$-divisor for short,
to be any coherent
sub ${\bfo}_{C}$-module  of $K^r=K^{\oplus r}$
having   rank $r$ and degree $n$.   Since $C$ is smooth, these
submodules are locally free and
coincide with the matrix divisors defined by A. Weil [31, 13].

Denote by $\div$ the set of all $(r,n)$-divisors.
This set can be identified with the set of rational points
of an algebraic ind-variety    $\divb$ that may be described  as follows.

For any  effective ordinary divisor $D$   set
$$
\div(D)=\{E\in \div\vert E\subseteq\bfo_{C} (D)^r\}
$$
where $\bfo_{C} (D)^r$ is considered
as a sub ${\bfo}_{C}$-module of $K^r$.
Clearly
$$
\div=\bigcup_{D\geq0} \div(D)
$$
and the elements of the set $\div(D)$ can be identified with the rational
points of the scheme  $\text{Quot}^{m}_{\bfo_{C} (D)^r/X/k}$,
$m=r\cdot\deg D-n$,
parametrizing torsion quotients
of $\bfo_{C} (D)^r$ having degree $m$ [14]. These
are smooth projective varieties and for any pair of effective divisors
$D,D'$ with $D\leq D'$ the inclusion
$$
\div(D)\rightarrow \text{Div}^{r,n}_{C/k}(D')
$$
is induced by a closed immersion of the corresponding varieties.

It is natural to stratify
the ind-variety
$\divb$ according to
Harder-Narasimhan type
$$
\divb = (\divb )^{ss} \cup \bigcup_{P\not= ss}\ \bfS_P
$$
where $(\divb )^{ss}$ is the open ind-subvariety of semistable divisors.
The cohomology of each stratum stabilizes and
this stratification is perfect.
( Here cohomology means $\ell$-adic cohomology for a suitable prime $\ell$. )
In particular, there is an identity
of Poincar\'e series
$$
P(\divb;t)= P((\text{\bf Div}^{r,n}_{C/k})^{ss};t) \ + \
\sum_{P\not = ss} P(\bfS_P;t)\cdot t^{2\d_P}\ .
\tag 1.1
$$
where $\d_P$ is the codimension of $\bfS_P$
( see Proposition 5.2 below for an explicit expression.~)

Let $r$ and $n$  be coprime. Then
the notions of stable
and semistable bundle over $C$ coincide, and
the   moduli space $N(r,n)$ of stable vector bundles   having
rank $r$ and degree $n$  is in this case a smooth
projective algebraic variety.
It is natural to define, by analogy with the
classical case, Abel-Jacobi maps
$$
\pmb\vartheta:(\text{\bf Div}^{r,n}_{C/k})^{ss}\longto N(r,n)\ ,
$$
by assigning to a divisor $E$ its isomorphism class as a
vector bundle.
There is a ``coherent locally free'' module $\pmb E$ over
$N(r,n)$, considered as a constant ind-variety, and a
morphism
$$
\matrix
({\text{\bf Div}}^{r,n}_{C/k})^{ss}& & @>\text{\bf i}>> & &\Bbb P(\pmb E )\\
_{\displaystyle\pmb\vartheta} &\searrow & & \swarrow & \\
&& N(r,n) &&
\endmatrix
$$
inducing an isomorphism in cohomology. It follows that
$$
H^*((\text{\bf Div}^{r,n}_{C/k})^{ss})=H^*(N(r,n))[x]
$$
where  $x$  is an independent variable of degree $2$, and hence
$$
P(N(r,n);t)= (1-t^2)\cdot P((\text{\bf Div}^{r,n}_{C/k})^{ss};t)\ .\tag 1.2
$$

If $(r'_1,d'_1),\dots,(r'_l,d'_l)$ is the type determined by $P$, then
the map
$$
\pmb\delta:\ (\text{\bf Div}^{r'_1,d'_1}_{C/k})^{ss}\times\dots\times
(\text{\bf Div}^{r'_l,
d'_l}_{C/k})^{ss}\longto\bfS_P\ ,
$$
given by $(E_1,\dots,E_l)\mapsto E_1\oplus\dots\oplus E_l$,
induces an isomorphism in cohomology and
one has
$$
P(\bfS_P;t)=\prod_{1\leq j\leq l}P(
(\text{\bf Div}^{r'_j,
d'_j}_{C/k})^{ss};t)\ .
\tag 1.3
$$

In order to find the Betti numbers of $N(r,n)$ it suffices, by (1.2), to
know those of $(\text{\bf Div}^{r,n}_{C/k})^{ss}$. But
the formulae (1.1) and (1.3)  reduce this computation,
for arbitrary $r$ and $n$, to the  calculation of those of
$\text{\bf Div}^{r,n}_{C/k} $.

The varieties $\div(D)$ are  analogous to Grassmannians and
share with  them the property of having a decomposition into
Schubert ``strata''  which may be defined, for example,
in terms of    a toric
action. It follows easily [3] that their cohomology
is free of torsion and that their Betti
numbers are given by
$$
P(\div(D);t)=\sum_{ \bold m} t^{2d_{\bold m}}\
P(C^{(m_1)};t) \cdots  P(C^{(m_r)};t)
$$
where $\bold m=(m_1, \ldots ,m_r)$ is any partition of $m=r\cdot\deg D-n$
by nonnegative integers, $d_{\bold m}=\sum_{1\le i\le r}(i-1)m_i$\ ,
and $C^{(m)}$ stands for the symmetric product.
This  reduces the calculation,  in a   sense, to the well-known  abelian
case $r=1$, and since  the cohomology of $\divb$ stabilizes,   one  obtains
$$
P(\divb;t)=\frac{\prod_{j=1}^r (1+t^{2j-1})^{2g}}
{ (1-t^{2r})
\prod_{j=1}^{r-1}(1-t^{2j})^2}\ .\tag 1.4
$$

There is a clear formal analogy between this approach and that of [2].
In [2] Atiyah and Bott consider the space $\Cal C(r,n)$ of
holomorphic structures on a fixed
$C^\infty$ vector bundle $E$ together with the action on it  of
the complexified gauge group
$\Cal G =\text{Aut }(E)$.
They stratify $\Cal C(r,n)$ according to Harder-Narasimhan type and,
using $\Cal G$-equivariant cohomology, obtain formul\ae, similar to
the ones above, that reduce the calculation of $P(N(r,n);t)$ to
that of $P(B\Cal G ;t)$ for arbitrary $r$. Since the Poincar\'e series
of the classifying space $B\Cal G$ coincides with the expression (1.4),
the two approaches give the same formul\ae\  for the Betti numbers as in
[24,15,17,9,30,18,6].

If one considers, for $k$ the
algebraic closure of a finite field $\Bbb F_q$, the
Abel-Jacobi map from the set of $\Bbb F_q$-rational points of $\divb$ to
the  (finite) set of isomorphism classes of $\Bbb F_q$-rational bundles,
having rank $r$ and degree $n$, then by simple considerations of
cardinalities one can prove the function-field analogue of the Siegel
formula [11].
Recall that this formula, equivalent to the Tamagawa number of
$\text{SL}(r,K)$ being 1, is the keystone  for the computation of the
Betti numbers of $N(r,n)$ in the arithmetic approach, via the Weil
conjectures, of Harder, Narasimhan, Desale and Ramanan [15,17,9].

\head 2. Some notation and preliminary remarks \endhead

Let $k$ be an algebraically closed field of arbitrary characteristic.
In this paper, an ind-variety
$\X=\{X_\lambda,f_{\lambda\mu}\}_{\lambda,\mu\in \Lambda}$
is simply an inductive
system of  $k$-algebraic varieties
indexed by some filtered ordered set $\Lambda$.
The indexing set $\Lambda$ will often be the set of effective
divisors on the curve $C$.
We shall say that  $\X $
is smooth if there is a $\lambda_0\in\Lambda$ such that
for every $\lambda\geq\lambda_0$, $X_\lambda$ is smooth.

Given another ind-variety
$\Y=\{Y_\rho,\  g_{\rho\sigma}\}_{\rho,\sigma\in\Gamma}$\ , a morphism
$\pmb\Phi=\{\alpha,\ \{\phi_\lambda\}_{\lambda\in\Lambda}\}$
from $\X$ to $\Y$ consists of an order preserving map
$\alpha :\Lambda\longto\Gamma$ together with a family of morphisms
$$
\phi_\lambda:X_\lambda\to Y_{\alpha(\lambda)}
$$
satisfying the obvious commutativity properties.
We shall say that a morphism $\pmb\Phi$ is a
quasi-isomorphism if:
\roster
\item  $\alpha(\Lambda)$ is a cofinal subset of
$\Gamma$;
\item for every integer $k$ there exists a $\lambda_k\in\Lambda$
such that for every $\lambda\geq\lambda_k$ the map $\phi_\lambda:X_\lambda
\to Y_{\alpha(\lambda)}$ is an open immersion and the codimension
of $Y_{\alpha(\lambda)}-\phi_\lambda(X_\lambda)$ in $Y_{\alpha(\lambda)}$
is greater than $k$.
\endroster

\sskip

Unless otherwise stated cohomology will mean $\ell$-adic
cohomology
$$
H^*(X)=\varprojlim_r H^*(X_{\text{\'et}}, \Bbb Z /\ell^r\Bbb Z )
$$
for a suitable prime $\ell$ distinct from the characteristic of the
field and coprime with the ranks  and the degrees   of the vector bundles
considered. If the base field is $\Bbb C$, then cohomology
could also be singular cohomology with coefficients in $\Z$ or in any
field ( with similar restrictions on the characteristic. )

Whenever we mention  the K\"unneth formula,
we are referring to  the one, without supports, in [8].

We  say that the cohomology of the ind-variety  $\X$
stabilizes if for every integer
$n$ there exists a $\lambda_n$
such that for every $\lambda,\mu$ with $\mu\geq\lambda\geq\lambda_n$ and
for every $j$ with $j\leq n$ the map
$$(f_{\lambda\mu})^*: H^j(X_\mu)\longto H^j(X_\lambda)$$
is an isomorphism.
The cohomology ring of
$\X $ will be
the projective limit
$$
H^*(X)=\varprojlim_{\lambda}
H^*(X_\lambda)\ .
$$
We  may also refer to the
ring $H^*(\X)$ as   the stable cohomology ring of $\X$.

The Poincar\'e series of $\X$ is by definition that
of $H^*(\X)$.

Note that if $\pmb\Phi:\X\to\Y$ is a quasi-isomorphism
of smooth ind-varieties, then  $\pmb\Phi^* :H^*(\Y)\to H^*(\X)$ is an
isomorphism ( see, for example,  Lemma 9.1 in Chapter VI of [23].~)
Furthermore, the cohomology of $\X$ stabilizes if, and
only if, that of $\Y$ does.

An $\bfo_{\X}$-module will be, in this paper, a family
$\pmb E =\{ \Cal E_{\lambda},\ \eta_{\lambda,\mu}\}_{\lambda ,\mu\in\Lambda}$
where $\Cal E_{\lambda}$
is a  $\bfo_{X_{\lambda}}$-module and
$$
\eta_{\lambda,\mu}:f^*_{\lambda,\mu}\Cal E_{\mu}\longto \Cal E_{\lambda}
$$
are surjective homomorphisms ( for example isomorphisms. )
To such a module, with coherent $\Cal E_{\lambda}$'s say,
one can associate ( as in E.G.A. )  an ind-scheme
$$
\Bbb P (\pmb E )\longto \X
$$
and a fundamental invertible module $\bfo_{\Bbb P} (1)$ over it.

Two special situations  will often arise.
The first one, is that
of a ``coherent locally free''
module of ``infinite rank'' over a ``finite dimensional'' base
$\X$. Here ``coherent locally free'' means that the coherent modules
$\Cal E_{\lambda}$ are locally free for $\lambda$ in a cofinal subset,
``infinite rank'' that the ranks of the
$\Cal E_{\lambda}$ are unbounded, and ``finite dimensional'' that the
cohomology of $\X$ stabilizes and, furthermore, vanishes beyond
a certain dimension. In this case it is easy to show that
$$
H^*(\Bbb P (\pmb E ) )=H^*(\X )[\xi ]
$$
where $\xi$, the Chern class of $\bfo_{\Bbb P} (1)$, is algebraically
independent.
The second
one, is that
of a ``coherent locally free'' module of ``finite rank''
( i.e. its components have bounded ranks.~)
In this case, there are Chern classes $c_i(\pmb E )\in
H^{2i}(\X)$, for $0\leq i\leq l$ say.
Note that $\bfo_{\Bbb P} (1)$ is already an example of this sort.

\head 3. The stratification associated to a direct sum decomposition \endhead

Let $X$ be a projective variety over a field $k$, and let $\Cal F$
be any coherent $\bfo_{X}$-module. Recall that there is a universal
family of coherent quotients of $\Cal F$ parametrized by a scheme
$Q=\text{Quot}_{\Cal F/X/k}$ which is a union of projective schemes
[14]. Since there is an obvious bijection between coherent
quotients of $\Cal F$ and coherent submodules of $\Cal F$, we can
also consider $Q$ as the variety parametrizing the universal family
of coherent submodules of $\Cal F$. We shall do this systematically
throughout the paper.

Recall, that if $E\subseteq\Cal F$  corresponds to a point of $Q$ then the
Zariski tangent space to $Q$ at $E$ can be identified with $\text{Hom}
(E,\Cal F/E)$. If   $\text{Ext}^1(E,\Cal F/E)=0$, then $Q$ is  smooth
at  $E$.

It is  often useful to think of the schemes
$\text{Quot}_{\Cal F/X/S}$ as a generalization of Grassmannians
( which correspond to the case $X=S$.~) One of the key techniques for studying
the geometry of Grassmannians is their decomposition into Schubert cells.
The classical way of defining these cells is via incidence conditions. There
is, however, a different approach via flows under the action of a suitable
torus
which yields the same results. From this point of view, one should
think of the results
in [3] as providing generalized Schubert ``strata'' in terms of a torus action.

Note that associated to any direct sum decomposition
$$
\Cal F=\Cal F_1\oplus\cdots\oplus\Cal F_n
$$
of a coherent module $\Cal F$, there
is a natural action of $\Bbb G_m^n$ on $\text{Quot}_{\Cal F/X/k}$ having
$$
\text{Quot}_{\Cal F_1/X/k}\times_k\cdots\times_k\text{Quot}_{\Cal F_n/X/k}
$$
as fixed points (proof as in [3].~)

For simplicity, we shall restrict ourselves to the case $X=C$ and $n=2$,
but one could develop the general case similarly along the lines of [3]
( see a paper, in preparation, by one of the authors. )

Choose  $\lambda_1,\lambda_2\in\Bbb Z$, with $\lambda_1 > \lambda_2$, and
define
a one-parameter subgroup
$\lambda:\Bbb G_m\to\Bbb G_m^2$ by
$\lambda (t)=(t^{\lambda_1},t^{\lambda_2})$. The induced $\Bbb G_m$-action
has the same fixed points as the original one, and
we now have a
decomposition of  $\text{Quot}_{\Cal F/X/k}$ into strata,
determined by the flow towards the different
components of the fixed points.

It is possible, however, to describe these strata in terms closer
to the classical description as follows.

Let $\Cal F_1, \Cal F_2$ be coherent $\bfo_C$-modules on the
projective curve $C$. A coherent submodule $E$  of $\Cal F=\Cal F_1
\oplus\Cal F_2$ determines submodules $E_1$ and
$E_2$ of $\Cal F_1$ and $\Cal F_2$, respectively, which make exact
the following commutative diagram
$$
\matrix
0&\to&\Cal F_1&\to&\Cal F_1\oplus\Cal F_2&\to&\Cal F_2&\to& 0\cr
&&\cup&&\cup&&\cup&&\cr
0&\to &E_1&\to &E&\to &E_2&\to& 0
\endmatrix
$$
Here, the horizontal maps are the obvious ones.

Furthermore, the map $E\subset\Cal F_1\oplus\Cal F_2\to\Cal F_1$
induces, after taking quotients with $E_1$, a homomorphism
$\varphi :E_2\to\Cal F_1/E_1$. Conversely,
given $E_1\subseteq \Cal F_1$, $E_2\subseteq\Cal F_2$
and $\varphi:E_2\to\Cal F_1/E_1$, we can recover $E$ as the fibre
product $\Cal F_1 \times _{\Cal F_1/E_1}E_2$.

Let $Q=Q^{(r,n)}$ be the scheme parametrizing
submodules of $\Cal F=\Cal F_1\oplus\Cal F_2$ having rank $r$ and
degree $n$. Consider, again, the action of $\Bbb G_m$ on $Q$ with the choices
made above.
Note that
the associated modules $E_1$,
$E_2$ do not change as $E$ describes a $\Bbb G_m$-orbit in $Q$.
Note also that,
while this happens,
$\varphi$ is being multiplied by $t^{\lambda_1 -\lambda_2}$.
It is now clear that the fixed points under this
action are the submodules of  type
$E=E_1\oplus E_2$ with $E_1\subseteq\Cal F_1$ and
$E_2\subseteq\Cal F_2$, and that one of the two
stratifications of $Q$ determined by the $\Bbb G_m$-action can
be described as
$$
Q^{(r,n)}=\bigcup\Sb r_1+r_2=r\\n_1+n_2=n\endSb\  Q^{(r_1,n_1;r_2,n_2)}
$$
where $Q^{(r_1,n_1;r_2,n_2)}$ consists of those $E$ such that the
associated modules $E_1$, $E_2$ have ranks $r_1$, $r_2$ and
degrees $n_1$, $n_2$ respectively.
Note that the other stratification is obtained after permutation of $\Cal F_1$
with
$\Cal F_2$.

For $i=1,2$\ , let $Q_i^{(r_i, n_i)}$ be the scheme parametrizing
$(r_i,n_i)$-submodules of $\Cal F_i$. There is
an obvious map
$$
Q^{(r_1, n_1;r_2,n_2)} \to Q_1^{(r_1, n_1)} \times Q_2^{(r_2, n_2)}
$$
whose fibre over $(E_1,E_2)$ can be identified with the vector space
$\text{Hom } (E_2,\Cal F_1/E_1)$.
If the dimension of the fibres is constant ( e.g. if
$\text{Ext}^1(E_2,\Cal F_1/E_1)=0$ for every $(E_1,E_2)\in
Q_1^{(r_1, n_1)}\times Q_2^{(r_2, n_2)}$ ), then $Q^{(r_1,n_1;r_2,n_2)}$
can be considered
as a vector bundle over  $Q_1^{(r_1, n_1)} \times Q_2^{(r_2, n_2)}$.

We shall refer to the stratification just described as being associated to the
direct sum decomposition
$\Cal F=\Cal F_1\oplus\Cal F_2$.

\head 4. The ind-variety  of divisors and its cohomology \endhead

Let $C$ be a smooth projective curve of genus $g$, defined
over an algebraically
closed field $k$, and let $\bfo_{C}$ be its structure sheaf.
We denote by   $K$
 the    constant
$\bfo_{C}$-module determined by the  function field of $C$.
\medskip

\definition{Definition 4.1}
A {\it divisor} of rank $r$ and degree $n$ over $C$, an $(r,n)$-divisor for
short,
is a  coherent  sub $\bfo_{C}$-module of $K^r$  having
degree $n$ and rank $r$.
We denote the set of all such divisors by $\div$.
\enddefinition

Note that, since $C$ is smooth,  coherent submodules of $K^r$
are locally free and these
divisors
coincide with the matrix divisors defined by A. Weil [31, 13].

For any effective ordinary divisor $D$  we have
$\bfo^r_{C}\subseteq \bfo_{C}(D)^r\subseteq K^r$ and  we  define
$$
\div(D)=\{E\vert E\in\div\text{ and }E\subseteq \bfo_{C}(D)^r)\}\ .
$$
On  each $\div(D)$ there are natural structures of smooth projective
algebraic variety.
In fact, the elements of   $\div(D)$ can be identified with the rational
points of the scheme  $\text{Quot}^{m}_{\bfo_{C} (D)^r/X/k}$,
$m=r\cdot\deg D-n$,
parametrizing torsion quotients
of $\bfo_{C} (D)^r$ having degree $m$ [14]. These
are smooth projective varieties and for any pair of effective divisors
$D$, $D'$ with $D\leq D'$ the inclusion
$$
\div(D)\rightarrow \text{Div}^{r,n}_{C/k}(D')
$$
is induced by a closed immersion of the corresponding varieties.\medskip

The ind-variety  $\divb$ is determined  by the
inductive system consisting of the varieties $\div(D)$ and   the
closed immersions above.

\remark{Remark}
One could also define, as in [4], an ind-variety  $\bold Q^{r,n}$ given by
$$
Q^{r,n}(D)=\text{Quot}^{n+r\cdot\deg D}_{\bfo_{C}^r/X/k}
$$
for every ordinary effective divisor $D$, and having as structure maps
the morphisms
$$
Q^{r,n}(D_1) \longto  Q^{r,n}(D_2)
$$
obtained by tensoring the submodules with $\bfo_{C} (-D)$, where
$D=D_2-D_1\geq 0$. It is easy to show that there is a natural isomorphism
of inductive systems
$$
\divb\longto \bold Q^{r,n} \ ,
$$
defined by  tensoring the submodules with $\bfo_{C} (-D)$ for every
ordinary effective divisor $D$.
\endremark

We are now ready to state the following

\proclaim{Proposition 4.2} The cohomology of $\text{\bf Div}^{r,n}_{C/k}$
is free of torsion, stabilizes and its Poincar\'e series is given by
$$
P(\text{\bf Div}^{r,n}_{C/k};t)=\frac{\prod_{j=1}^r (1+t^{2j-1})^{2g}}
{ (1-t^{2r})
\prod_{j=1}^{r-1}(1-t^{2j})^2}\ .
$$
\endproclaim

\demo{Proof}
We shall prove below that the cohomology stabilizes.

Recall ([3]; cf. also the remark at the end
of this section) that the cohomology  of $\div(D)$
is free of torsion and that its Poincar\'e series is given by
$$
P(\div(D);t)=\sum_{ \bold m} t^{2d_{\bold m}}\
P(C^{(m_1)};t) \cdots  P(C^{(m_r)};t)
$$
where $\bold m=(m_1, \ldots ,m_r)$ is any partition of
$ m  = r\cdot\deg D-n$
by non-negative integers and $d_{\bold m}=\sum_{1\le i\le r}(i-1)m_i$.

Consider the formal power series $E(t,u)$ defined by
$$
E(t,u)=\sum_{m\geq0}\left(\sum_{ |\bold m |=m}P(C^{(m_1)},t)
\cdots  P(C^{(m_r)},t)\cdot t^{2 d_{\bold m}} \right)\cdot u^m
$$
where $|\bold m |= \sum m_i$.

We know, see [21],  that
$$
\sum_{j\geq0}P(C^{(j)},t)u^j={(1+ut)^{2g}\over (1-u)(1-ut^2)}\ .
$$
Therefore
$$E(t, u)=\prod_{j=0}^{r-1} {(1+u\ t^{2j+1})^{2g}\over
(1-u\ t^{2j})(1-u\ t^{2j+2})}\ .
$$
is a rational function, and
$$
P(\divb;t)=-\operatornamewithlimits{Res}_{u=1}E(t,u)=
\frac{\prod_{j=1}^r (1+t^{2j-1})^{2g}}
{ (1-t^{2r})
\prod_{j=1}^{r-1}(1-t^{2j})^2}\
$$
has the expected value.

Let us now prove the claim, about stability of the cohomology, made
at the beginning.

\proclaim{Proposition 4.3} If $D'\geq D$, then the morphism
$$
H^i(\div (D')\longto H^i(\div(D))
$$
induced by the closed immersion $\div(D)\rightarrow \div(D')$
is an isomorphism for $0\leq i<r\deg D-n$.
\endproclaim

This proposition is an immediate consequence of the next lemma,
since we can always factor  the map $\div(D)\rightarrow \div(D')$ through
a succession of maps of the type considered there.

\proclaim{Lemma 4.4} Let $\L$, $\F$ be  locally free $\bfo_{C}$-modules
having ranks $1$  and $r-1$ respectively.  Fix a point $P$ in $C$ and let
$Q$, $Q'$ be the varieties of $(r,n)$-submodules
of $\L\oplus \F$  and   $\L(P)\oplus\F$ respectively. Then the morphism
$$
H^i(Q')\longto H^i(Q),
$$
induced by the closed immersion $Q\rightarrow Q'$, is an isomorphism
for $0\leq i<\deg\F+\deg\L-n$.
\endproclaim

\demo{Proof} Suppose first that $r=1$. We can take $\L=\bfo_{C}
(D)$ for some divisor $D$. Every submodule of $\bfo_{C}(D)$ is equal to
$\bfo_{C}(F)$ for some  $F\leq D$. The map \ $\bfo_{C} (F)\mapsto D-F$\  gives
an
isomorphism \  $Q\longto C^{(m)}$, where $m=\deg\L-n$.
Similarly, there is an isomorphism \ $Q'\longrightarrow
C^{(m+1)}$ \ and a commutative diagram
$$
\CD
Q@>>> C^{(m)}\\
@V VV @V VV\\
Q' @>>> C^{(m+1)}
\endCD
$$
where the map \ $C^{(m)}\rightarrow C^{(m+1)}$ \
is given by $H\mapsto H+P$. As is well known [21], the induced   map
$$
H^i(C^{(m+1)})\longto H^i(C^{(m)})
$$
is an isomorphism for $0\leq i< m$,
and a monomorphism for $i= m$.
Therefore, the same holds for  the map $H^i(Q')\longto H^i(Q)$.
\enddemo

In general,   consider the stratification of $Q$
associated to the  direct sum
$ \L\oplus\F  $.
We can write
$$
Q=Q_0\cup Q_1\cup Q_2\cup\dots
$$
where $Q_j$ consists of the submodules $E$ of $\L\oplus\F$ which
project onto a submodule $E_2$ of $\F$ having degree $\deg\F-j$.
The subvariety $Q_j$ has codimension $j$ in $Q$ and can be considered
as a vector bundle over $V_j\times Z_j$, where $V_j$ is the
variety of submodules of $\L$ having rank 1 and degree $n-\deg \F+j$,
and $Z_j$ is the variety of submodules of $\F$ having rank $r-1$
and degree $\deg\F-j$.

Analogously, we have a stratification
$$
Q'=Q'_0\cup Q'_1 \cup Q'_2\cup\dots
$$
The subvariety $Q'_j$ is now a vector bundle over $V'_j\times Z_j$, where
$Z_j$ is as above and
$V'_j$ is the variety of submodules of $\L(P)$ having rank 1 and degree
$n-\deg \F+j$.

We have, for every $j$,  a commutative diagram of Gysin exact
sequences
$$
\CD
\to  H^{p-2j}(Q'_j) @>>> H^p(Q'_0\cup\dots\cup Q'_j) @>>>
H^p(Q'_0\cup\dots\cup Q'_{j-1})\to \\
@VVV @VVV @VVV  \\
\to H^{p-2j}(Q_j) @>>> H^p(Q_0\cup\dots\cup Q_j)
@>>> H^p(Q_0\cup\dots
\cup Q_{j-1})\to
\endCD
$$
We know that the map $H^i(V'_j)\to H^i(V_j)$ is an isomorphism
for $0\leq i<\deg\L+\deg\F-n-j$ and a monomorphism for $i=\deg
\L+\deg\F-n-j$. Using the isomorphisms
$$
H^*(Q'_j)\simeq H^*(V'_j)\otimes H^*(Z_j)
$$
and
$$
H^*(Q_j)\simeq H^*(V_j)\otimes H^*(Z_j)
$$
we see that
$H^i(Q'_j)\to H^i(Q_j)$ is an isomorphism for $0\leq i<\deg
\L+\deg \F-n-j$ and a monomorphism for $i=\deg\L+\deg\F-n-j$.

At this point the five-lemma applied to the diagram above allows us
to perform an induction on $j$  and conclude the proof of the lemma.
\enddemo
\remark{Remark}. Since the stratification $Q=Q_0
\cup Q_1\cup\dots$ comes
from an action of $\Bbb G_m$, it is easy to show that it is perfect
i.e. the bottom
Gysin sequence appearing in the diagram above
breaks up into short exact sequences, one
for every $j$.
It is  now possible, by induction on $r$, to show
that the cohomology of $\div(D)$ is  free of torsion, and to compute
its Poincar\'e polynomial.
\endremark

\head 5. The Shatz stratification of the ind-variety of divisors \endhead

Let us recall some standard definitions. The slope $\mu(E)$
of a coherent $\bfo_{C}$-module $E$ is defined by
$$
\mu(E)={ \deg E \over\rank E}\ .
$$
A coherent locally free $\bfo_{C}$-module $E$ is called
semistable (resp. stable) if  for
every proper submodule $F$ of $E$ one has
$$
\mu(F)\leq \mu(E)\qquad\qquad \text{(resp. } \mu(F)<\mu(E))\ .
$$
Recall the following two properties of semistable modules.
\roster
\item"(5.1)" If $E$ is semistable and $\mu(E)>2g-1$, then $E$ is
generated by its global sections and $H^1(C,E)=0$.
\item"(5.2)" If $E$ and $F$ are semistable, then
$\text{Hom }(E,F)=0$ whenever  $\mu(E)>\mu(F)$.
\endroster

The Harder-Narasimhan filtration of a coherent locally free $\bfo_{C}$-module
$E$ is the unique filtration
$$
0=E_0\subset E_1\subset\dots\subset E_l=E
$$
with $G_i=E_i/E_{i-1}$  semistable and
$$
\mu(G_1)> \mu(G_2)>\dots >\mu(G_l) \ .
$$
The type  of $E$ will be the sequence of pairs
$\{ (\rank G_i,\deg G_i)\}_{1\leq
i\leq l}$. Equivalently, the Shatz polygon of $E$ will have  the points
$\{ (\rank E_i,\deg E_i)\}_{0\leq i\leq l}$ as vertices. Note that if
$\{ (r_i,d_i) \}_{0\leq i\leq l}$ is the polygon of $E$,
then the type of $E$
is given by
$\{ (r'_i,d'_i)\}_{1\leq i\leq l}$ \ with
$$
(r'_i,d'_i)=(r_i-r_{i-1},d_i-d_{i-1}) .
$$
We shall
use primed letters for the type of a bundle, and unprimed letters for
the vertices of its polygon.

Recall that the Shatz  polygon  $P_E$ of  $E$ is characterized by the
following property.
\roster
\item"(5.3)" If $F$ is a $(s,m)$-submodule of $E$,
then the point $(s,m)$ lies  either on  $P_E$ or below it.
\endroster

Consider the set $\Cal P_{r,n}$  of all strictly convex polygons
of $\Bbb R^2$ joining $(0,0)$ to $(r,n)$. If $P,$ and $P'$ are in
$\Cal P_{r,n}$, then we say that $P\geq P'$ whenever $P$ lies above $P'$.

Let  $\Cal E$ be a family of vector bundles of rank $r$ and degree
$n$ over $C$ parametrized by a scheme $T$. For $t\in T$ let $P_t$
be the polygon of $\Cal E_t$. For each $P\in \Cal P_{r,n}$  we define
subsets of $T$
$$
F_P(T)=\{t\in T\vert  P_t>P\}\quad , \quad \Omega_P(T)=T-F_P(T)
$$
and
$$
S_P(T)=\{t\in T\vert P_t=P\}\ .
$$
We may also write $(T)^P$, instead of $S_P(T)$.

The partition
$$
T=\bigcup_{P\in\Cal P_{r,n}} S_P(T)
$$
is called the Shatz stratification of $T$.

We summarize, in the next proposition, some general properties of this
stratification in the case $T=\div(D)$. In order to state  them we shall
need some notation.

\definition{Definition 5.1}
Let $X$ be a scheme and let
$$
0=\F_0\subseteq\F_1\dots\subseteq\F_l=\F
$$
be a filtration
of an $\bfo_{X}$-module $\Cal F$. We
denote by $\Cal H om_-(\F,\F)$ the subsheaf of $\Cal H om
(\F,\F)$ consisting of the germs of homomorphism  preserving
the filtration. We also set
$$
\Cal H om_+(\F,\F)={\Cal H om(\F,\F)/\Cal H om_-(\F,\F)} \ .
$$
\enddefinition

\proclaim{Proposition 5.2}
Let $T_D=\div(D)$ and
fix a  Shatz polynomial
$$
P=\{ (r_0,d_0),\dots ,(r_l,d_l)\} \ .
$$
We set, as usual, $r'_i=r_i-r_{i-1}$ and $d'_i=d_i-d_{i-1}$. Then
\roster

\item"(1)"  $F_P(T_D)$ is closed and  $\Omega_P(T_D)$ is an open set;

\item"(2)"  $S_P(D)=S_P(T_D)$ is in a natural way a closed subscheme of
$\Omega_P
(T_D)$. Moreover, if $\Cal U^P_D$ is the restriction
to $S_P(D)$ of the universal family of divisors,
then $\Cal U^P_D$ has a universal Harder-Narasimhan
filtration
$$
0=\Cal U^P_0\subset\Cal U^P_1\subset\dots\subset\Cal U^P_l=\Cal U^P_D
$$
such that
$$
\Cal U^P_i/ \Cal U^P_{i-1}
$$
is locally free for  $i=1,\dots, l$.

\item"(3)"  If $\deg D>{\displaystyle d'_1\over\displaystyle r'_1}+2g-1$,
then $S_P(D)$ is smooth,
has codimension
$$
d_P=\sum_{i>j}
r'_ir'_j\left[\left({\displaystyle d'_j\over\displaystyle r'_j}-
{\displaystyle d'_i\over\displaystyle r'_i}\right)+g-1\right]
$$
and its normal bundle
can be identified with
$$
R^1q_*\ \Cal H om_+(\U^P_D,\ \U^P_D )
$$
where $q:C\times S_P(D)\longto S_P(D)$ is the
natural projection.

\item"(4)" If $\deg D\leq{\displaystyle d'_1\over\displaystyle r'_1}+2g-1$,
then
$$
\text{\rm codim }S_P(D)\geq\deg D-c
$$
where $c$ is a constant independent of $P$ and $D$.
\endroster
\endproclaim

\demo{Proof} Most of the assertions made in the proposition are easy
consequences of general results [29,30];
however, for the convenience of the reader, we outline a
proof.

Note that throughout this outline, if $E$ is a submodule of $\bfo_{C}(D)^r$,
we shall often write $\EE$
instead of $\bfo_{C}(D)^r/E$.

Let $\Cal D^{s,m} _D=\Cal D^{s,m}_r(D)$ be the projective scheme parametrizing
the
$(s,m)$-submodules of $\bfo_{C}(D)^r$.
We shall show that the proposition holds, more generally, for these varieties.

Let $P$ be any polygon with vertices
$$
(0,0)=(s_0,m_0), \dots ,(s_l,m_l)=(s,m) \ .
$$
We denote by $\text{Drap}^P_D$ the projective scheme parametrizing
the universal family of flags
$$
0=E_0\subset E_1\dots\subset E_l\subset\bfo_{C}(D)^r
\tag 5.4
$$
having $P$ as associated polygon i.e.
$$
(\rank E_i, \deg E_i)= (s_i,m_i)
$$
for $i=0,\dots ,l$.
The existence of $\text{Drap}^P_D$ follows, by induction on $l$,
from that of the schemes of quotients.

In fact $\text{Drap}^P_D$ can be identified, in a natural
way, with a closed subscheme of
$$
\Cal D^{s_1,m_1}_D\times\dots\D^{s_l,m_l}_D \ .
$$
The Zariski tangent space to $\text{Drap}^P_D$ at the point
corresponding to a flag like (5.4) can be identified
with the subspace of
$$
\text{Hom }(E_1,\EE_1)\oplus\dots\oplus\text{ Hom }
(E_l,\EE_l)
$$
consisting of the $l$-tuples $(f_1,\dots , f_l)$ such that
$f_{i+1}$ and $f_i$ map to the same element in $\text{Hom }
(E_i,\EE_{i+1})$ under the obvious maps.

Let us begin the proof of $(1)$.
Consider the map
$$
\pi_P:\text{Drap}^P_D\longto\Cal D^{s,m}_D
$$
which assigns to a filtration, as in (5.4), the submodule $E_l$.
Since $\pi_P$ is projective, its image is closed in $\Cal D^{s,m}_D$.
Also, as a consequence of (5.3), we have for every $(s,m)$-polygon $P$
$$
F_P(\Cal D^{s,m}_D)=\bigcup_{P'>P}
\text{Im}\,\pi_{P'}\ .
$$
Since the degree of any submodule of $\bfo_{C}(D)^r$ is bounded by $r\deg D$,
there are only finitely many polygons $P'>P$ with nonempty
$\text{Drap}^{P'}_D$. This proves $(1)$.

Let now $U=\Omega_P(\Cal D^{s,m}_D )$ and consider
$$
\pi'_P:\pi^{-1}_P(U)\longto U \ ,
$$
the restriction of $\pi_P$. This map is still projective,
and its image coincides with $S_P(\D^{s,m}_D)$ and is closed in
$U$. The uniqueness of the Harder-Narasimhan filtration shows that
$\pi'_P$ is injective. We claim that it is actually a closed immersion.
What needs to be shown ( cf. [30] ) is that if $(f_1,
\dots f_l)$ belongs to the tangent space to $\text{Drap}^P_D$
at a filtration  like (5.4) which also happens to be the
Harder-Narasimhan filtration of $E_l$,
then $f_i$ is uniquely determined by $f_{i+1}$ for $i=1,\dots,
l-1$.
But note that,
as a consequence of (5.2), one has $\text{Hom }(E_i,E_{i+1}/E_i)=0$;
therefore the map
$$
\text{Hom }(E_i,\EE_i)\to\text{Hom }(E_i,\EE_{i+1})
$$
is injective and $f_{i+1}$ uniquely
determines $f_i$.

We can  now identify $S_P
(\D^{s,m}_D)$ with $\pi^{-1}_P(U)$, an open
subscheme of $\text{Drap}^P_D$, thus proving $(2)$.

Given a polygon $P$, let $P_i$, for $0\leq i\leq l$, be the polygon
having as vertices the first $i+1$ vertices $(r_0,d_0)$, $\dots$,
$(r_i,d_i)$ of $P$.
Consider the sequence of morphisms
$$
S_{P_l}(\D^{r_l,d_l}_D)\buildrel \phi_l\over\longto
\dots\to S_{P_i}(\D^{r_i,d_i}_D)\buildrel \phi_i\over\longto
 S_{P_{i-1}}(\D^{r_{i-1},d_{i-1}}_D)\to\dots\longto S_{P_0}(\D^{0,0}_D)
\tag 5.5
$$
defined by
$$
\phi_i(E_i)=E_{i-1}
$$
whenever
$$
0=E_0\subset E_1 \subset\dots\subset E_{i-1}\subset E_i
$$
is the Harder-Narasimhan filtration of $E_i$.

It is easily verified, using the description of the tangent spaces to the
strata given above, that the differential of each $\phi_i$ at $E_i$
$$
(\d\phi_i)_{E_i}:T_{E_i}\longto T_{E_{i-1}}
$$
fits into an exact sequence
$$
0\to\text{ Hom }(E_i/E_{i-1},\EE_i)\longto T_{E_i}
\longto T_{E_{i-1}}\longto Q (E_i,E_{i-1})\to 0
\tag 5.6
$$
where $Q(E_i,E_{i-1})\subset\text{Ext}^1(E_i/E_{i-1},
\EE_i)$ is the image of
the composite of the obvious maps
$$
\text{Hom }(E_{i-1}, \EE_{i-1})\longto\text{ Hom }
(E_{i-1},\EE_i)
$$
and
$$
\text{ Hom }(E_{i-1},\EE_i)
\longto\text{ Ext}^1(E_i/E_{i-1},\EE_i)\ .
$$
The vector space $\text{Hom }(E_i/E_{i-1},\EE_i)$ can be
identified with the tangent space to the fibre $\phi^{-1}_i(E_{i-1})$ at $E_i$.
Indeed, the fibre $\phi^{-1}_i(F)$ over any $F$ is isomorphic to an open
subscheme of the projective scheme parametrizing $(r'_i,d'_i)$-submodules
of $\widetilde F$ ( namely, the subscheme consisting
of the semistable ones. )

Make now the assumption that $\deg D>\dfrac{d'_1}{r'_1}+2g-1$. Since $P$ is
convex,
this is equivalent to the assumption that $\deg D>\dfrac{d'_i}{r'_i}+2g-1$ for
$i=1,\dots, l$.

Then, as a consequence of (5.1), we have
$$
\text{Ext}^1(E_i/E_{i-1}, \EE_i)=0\ .
$$
This implies,  that the non-empty
fibres of $\phi_i$ are all smooth and  have the same dimension. Moreover,
the exact sequence above ( where $Q(E_i,E_{i-1})$ now vanishes ) shows
that the differential of $\phi_i$ is everywhere surjective.
It is also not difficult to verify that $\phi_i$ is surjective
 ( observe that if
$E_{i-1}$ is a submodule of $\bfo_{C} (D)^r$ having rank $r_{i-1}$, then
there is an injection
$\bfo_{C} (D)^{r-r_{i-1}}\rightarrow \EE_{i-1}$
and one can use (5.1) to inject any semistable $(r'_i,d'_i)$-bundle
into $\bfo_{C}(D)^{r'_i})$.
It follows, by induction on $i$, that the strata $S_{P_i}(\D^{r_i,s_i}_D)$
are smooth and that each $\phi_i$ is smooth.
The equality $\text{codim }S_{P_i}(\D^{r_i,s_i}_D)=\d_{P_i}$
can also be verified by induction on $i$ using the
Riemann-Roch theorem and the exact sequence (5.6).

To complete the proof of $(3)$ we still have to verify  the assertion made
there about the normal bundle. We shall
limit ourselves here  to showing that the normal space
$N_E$  to $S_P(\D^{s,m}_D)$
at a point $E$   is naturally isomorphic to $H^1(\Cal H
om_+(E,E))$, where $\Cal H om_+(E,E)$ is taken  relative to the
Harder-Narasimhan  filtration of $E$.

We already know that $N_E$ is isomorphic to $\text{Hom}(E,\EE)/
T_E$, where $T_E$ is the subspace of $\text{Hom }(E,\EE)$ consisting
of the homorphisms $f$ such that for $i=1,\dots, l$ there are
homorphisms $f_i\in\text{Hom }(E_i,\EE_i)$, with
$f_l=f$, making the following diagrams commutative
$$
\CD
E_i @> f_i >> \EE_i\\
@V  VV @V  VV\\
E_{i+1} @>f_{i+1}>> \EE_{i+1}
\endCD
$$
for $i=1,\dots , l-1$.
Consider now the exact sequences
$$
H^1(\Cal H om_-(E,E))\buildrel \alpha\over\rightarrow H^1
(\Cal H om(E,E))\longto H^1(\Cal H om_+(E,E))\to 0
$$
and
$$
\text{Hom}(E,\EE)\buildrel\beta\over\longto \text{ Ext}^1
(E,E)\longto\text{Ext}^1(E,\EE)\ .
$$
The assumption that  $\deg D>\dfrac{d'_i}{r'_i}+2g-1$ together with (5.1) imply
$$
\text{Ext}^1(E,\EE)=0\ ,
$$
and  $\beta$ is surjective.
Therefore we will know that $\beta$ induces an isomorphism
$$
\text{Hom}(E,\EE)/T_E\buildrel\sim\over\longto
H^1(\Cal H om_+(E,E))
$$
as soon as we show that $\beta^{-1}(\text{Im}\,\alpha)=
T_E$.
But this is clear once we interpret $\text{Im}\,\alpha$ (resp. $T_E$)
as the space
of first order infinitesimal deformations of $E$ (resp.
of $E$ inside $\bfo_{C}(D)^r$) which,
at least in one way, deform the Harder-Narasimhan filtration of $E$.
(In fact there is at most one way to do this:
$(f_1,\dots,f_l)$ is, as we know, uniquely determined
by $f$ and, similarly, one can show that $\alpha$ is injective.~)
This completes the proof of $(3)$.

We claim that for every $E\in \D^{s,m}_D$, with Harder-Narasimhan filtration
as in (5.4), one has
$$
\sum^l_{i=1}\text{ dim }\text{Ext}^1(E_i/E_{i-1},\EE_i)\le
c\tag 5.7
$$
where $c$ is a constant that depends  only on $r$.
Observe that
$$
\dim\text{Ext}^1(E_i/E_{i-1},\EE_i)\le r\dim H^1
((E_i/E_{i-1})^\vee\otimes \bfo_{C}(D))\ .
\tag 5.8
$$
Also
$$
\eqalign{
\mu((E_i/E_{i-1})^\vee\otimes \bfo_{C} (D))&=\deg D-\dfrac{d'_i}{r'_i}\ge \cr
& \ge\deg
D-\dfrac{d_1}{r_1}\geq0\ .\cr }
$$
and, if
$$
H^1((E_i/E_{i-1})^\vee\otimes \bfo_{C} (D))\not=0
\tag 5.9
$$
then, because of (5.1),
one has
$$
\mu((E_i/E_{i-1})^\vee\otimes \bfo_{C} (D))\le 2g-1 .
$$
It follows that
the modules $(E_i/E_{i-1})^\vee\otimes \bfo_{C} (D)$
satisfying (5.9)
vary in a bounded family. This together with (5.8)
implies (5.7). Now, taking into account the exact sequence (5.6) and
using (5.7)  we  conclude that
$$
\text{codim }S_P(\D^{s,m}_D)\ge \d_P-c\ .\tag 5.10
$$
Because of the way $\d_P$ is defined and because $P$ is convex, we also have
$$
\d_P\ge \dfrac{d'_1}{r'_1}-\dfrac{d'_l}{r'_l}-1\ge
\dfrac{d'_1}{r'_1}-\dfrac{m}{s}-1\ .
$$
If  $\deg D\le \dfrac{d'_1}{r'_1}+2g-1$, then we obtain
$$
\d_P-c\ge \deg D-c-{m\over s}-2g\ .
$$
This inequality, together with (5.10), shows that $(4)$
holds  and we are done.
\enddemo
\head 6. A digression on moduli spaces  of rigidified semistable
bundles \endhead

Although there is not even a coarse moduli space for semistable bundles,
there is always a fine moduli space for  rigidified semistable bundles.
We shall need this result in the following section.

Let $\x$ be a finite nonempty set of
points on the curve $C$. A family, parametrized by a scheme $S$,  of
$\x$-rigidified vector
bundles of rank $r$ over $C$ is a pair $(E,u)$ where $E$ is a vector
bundle of rank $r$ over $C\times S$ and $u$ is an isomorphism
$\bfo_{X}^r\to E|_X$, $X=\x\times S$.

\proclaim{Proposition 6.1} There
exists a universal family of $\x$-rigidified
semistable bundles of rank $r$ and degree
$d$. We shall denote  $M(r,d\,;\, \x )$ the parameter space of
this  family.
\endproclaim

\demo{Proof} As is well known [25, 26] there is a family $\Cal F$
of semistable $(r,d)$-bundles over $C$ which is
parametrized
by a smooth quasi-projective variety $T$ and  has
compatible actions of $\text{GL}(N)$ on
both $\F$ and
$T$ having  the following properties:
\roster
\item"a)" $\F$ has the local universal property for families of semistable
bundles of rank $r$ and degree $d$ over $C$;
\item"b)" If $t_1,t_2\in T$, then $\F_{t_1}$ is isomorphic to $\F_{t_2}$
if, and only if, $t_1$ and $t_2$ are in the same orbit;
\item"c)" the stabilizer of $t\in T$ maps isomorphically onto the group
of automorphisms of $\F_t$;
\item"d)" a good quotient of $T$  under the action of $GL(N)$ exists.
\endroster

For each $x\in\x$, let $q_x:\Cal B(x)\longto T$
be the fibre bundle whose fibre over
$t\in T$ is naturally isomorphic to the variety of all bases of the
vector space $\F (x,t)$.
Set
$$
\Cal B(\x )={\prod_{x\in\x} }\,_{\tsize T}\ \Cal B(x)
$$
and
$\tilde{\F}=(1_C\times q)^*\F$, where $q:\Cal B(\x )\longto T$
is the natural projection.

Now $\tilde{\F}$ is canonically  $\x $-rigidified and the action
of $GL(N)$ on $\F$ and  $T$ lifts to an action on $\tilde{\F}$ and
$\B(\x )$.
Since   the map $q$ is affine and because of property
d),  a good quotient of $\B(\x )$  under
the action of $GL(N)$ exists.
( This will be $M(r,d \,;\, \x )$. )
An automorphism of a semistable bundle $E$  that leaves fixed a basis
of $E(x)$, for $x$ a point in $C$,  has to be  the identity.
It follows from this and from property c) that the stabilizer of
any $z\in\B (\x )$  is the identity. Therefore $M(r,d\, ;\, \x )$
is a geometric quotient.

Recall that if   $\L$ and $\N$ are  two families of semistable
$\x$-rigidified bundles parametrized
by the same $S$, then the set
$$
S'=\{s\in S \vert\text{ $\L_s$  is isomorphic to $\N_s$
as  $\x$-rigidified bundles } \}
$$
is closed
in $S$ and there is an isomorphism of families of  $\x$-rigidified
bundles
$$
\L\vert_{C\times S'}\buildrel\sim\over\longto \N\vert_{C\times S'}\ .
$$

It follows that the action of $GL(N)$ on $\B(\x )$
is free,   and hence that $\B(\x )$ is a $GL(N)$-principal bundle
over $M(r,d\,;\, \x )$. At this point,   standard
criteria for descent allow us to    conclude that $\tilde{\F}$
descends to a vector bundle $\E $ over $C\times M(r,d\,;\, \x )$.
Finally, using properties a) and b)  one  can verify
that $\E$ is indeed a universal family of $\x$-rigidified semistable
bundles of rank $r$ and degree $d$.
\enddemo
\head 7. Cohomology of the strata \endhead

Let $P$  be a Shatz polygon having vertices $(r_0, d_0), \ldots ,(r_l,d_l)$.
We set, as usual,
$r'_i=r_i-r_{i-1}$ and $d'_i=d_i-d_{i-1}$  for  $i=1,\dots, l$.
Consider the   closed
immersion of ind-varieties
$$
\pmb\delta:(\text{\bf Div}^{r'_1,d'_1}_{C/k})^{ss}\times\dots\times
(\text{\bf Div}^{r'_l, d'_l}_{C/k})^{ss}\longto
\bfS_P\ .
$$
defined by
$\pmb\delta ((F_1,\dots, F_l))=F_l\oplus F_{l-1}\oplus\dots\oplus F_1$.

The aim of this section is to prove the following

\proclaim{Proposition 7.1}  The cohomology of the stratum
$\bfS_P$ stabilizes.  Furthermore,  the morphism $\pmb\delta$ induces
an isomorphism in cohomology, and we have an identity of
Poincar\'e series
$$P(\bfS_P ;t)=\prod^l_{i=1}P((\text{\bf Div}^{r'_i,d'_i}_{C/k})^{ss};t)\ .$$
\endproclaim
Before we start proving this proposition it will be useful
to have a preliminary result.

\proclaim{Proposition 7.2}
Let $\D^{s,m}_r(D)$ be the projective scheme
parametrizing $(s,m)$-submodules of $\bfo_{C}(D)^r$.
Also, let $\C^{s,m}_r(D)$ be  its closed subset  consisting of
the submodules $L$ of $\bfo_{C}(D)^r=\bfo_{C}(D)^{r-s}\oplus
\bfo_{C}(D)^s$ which do not project injectively into $\bfo_{C}(D)^s$.

We have
$$
\dim \D^{s,m}_r(D)\leq r( s\deg D-m )+ a\ ,
\tag 7.1
$$
and
$$
\dim\C^{s,m}_r(D)\leq s(r-1)\deg D+ b
\tag 7.2
$$
where both $a$ and $b$  are  constants independent of $D$.

Moreover,   for  $P\not= ss$,
define $\bfS_P^o$  as the open subset of the
$P$-stratum consisting of   divisors $E$  having  Harder-Narasimhan
filtration
$$
0=E_0\subset\dots\subset E_{l-1}\subset E_l=E
$$
such that
$E_{l-1}$ projects injectively into the second factor
of the direct sum decomposition
$\bfo_{C}(D)^r=\bfo_{C}(D)^{r-r_{l-1}}\oplus\bfo_{C}(D)^{r_{l-1}}$.

Then the open immersion
$$
\bfS_P^o\rightarrow \bfS_P
$$
is a quasi isomorphism ( in particular
it induces an isomorphism in cohomology.~)
\endproclaim

\demo{Proof}  The Zariski tangent space to the variety $\D^{s,n}_r(D)$
at a point
$L$ is isomorphic to   $\text{Hom }(L, \bfo_{C}(D)^r/L)$.
We can bound
the dimension of this vector space from above using
the Riemann-Roch theorem
and the obvious inequalities
$$
\eqalign{
\dim\text{Ext}^1(L,\bfo_{C}(D)^r/L)&\leq\dim
\text{Ext}^1(L,\bfo_{C}(D)^r)\cr
&\leq
\dim\text{Ext}^1(\bfo_{C}(D)^r,\bfo_{C}(D)^r)\leq r^2g\ .\cr}
$$
The resulting bound gives  the first inequality.

Consider now the stratification of $Q=\D^{s,m}_r(D)$ associated
to the  direct sum decomposition
$$
\bfo_{C}(D)^{r}= \bfo_{C}(D)^{r-s}\oplus\bfo_{C}(D)^s
$$
The stratum \ $Q^{(s_1,m_1;s_2,m_2)}$, where    $s_1+s_2\!=\!s$ and
$m_1+m_2=m$,  maps to $\D^{s_1,m_1}_{r-s}(D)\times\D^{s_2,m_2}_s(D)$
and the fibre over $(L_1,L_2)$  is isomorphic to $\text{Hom}(L_1,\bfo_{C}
(D)^{r-s}/L_2)$.

We can bound the dimension of
$\D^{s_1,m_1}_{r-s}(D)\times \D^{s_2,m_2}_s(D)$
using  (7.1), and  that of  the vector space
$\text{Hom}(L_1,\bfo_{C}(D)^{r-s}/L_2)$
using both the Riemann-Roch theorem and the obvious inequalities
$$
\eqalign{
\dim\text{Ext}^1(L_2,\bfo_{C}(D)^{r-s}/L_1)&\leq\dim\text{Ext}^1
(L_2,\bfo_{C}(D)^{r-s})\leq\cr
&\leq\dim\text{Ext}^1(\bfo_{C}(D)^s,\bfo_{C}(D)^{r-s})\leq
s(r-s)g\ .\cr}
$$
As a result,  we obtain
$$
\dim Q^{(s_1,m_1\,;\,s_2,m_2)}\leq (r-s+s_2)s\deg D-(r-s+s_2)m+b\ ,
$$
where $b$ is a constant  that depends  only on $r$.
Since $\C^{s,n}_r(D)$ is, by definition, the union of the strata for which
$s_2\leq s-1$, the second inequality   follows.

Finally, recall the map
$$
\phi_l:S_P(D)\longto S_{P_{l-1}}(\D^{r_{l-1},d_{l-1}}_r(D))
$$
used in  the proof of Proposition 5.1
If $E\in S_P(D)-S_P^o(D)$, then by definition
$$
\phi_l(E)\in S_{P_{l-1}}(\D^{r_{l-1},d_{l-1}}_r(D))
\cap \Cal C^{r_{l-1},d_{l-1}}_r(D)\ .
$$
We know that, for $\deg D$ large enough,  $\phi_l$ is a smooth surjective map
between smooth varieties, and we have
$$
\dim S_{P_{l-1}}(\D^{r_{l-1},d_{l-1}}_r(D))=r_{l-1}r\deg D+c
$$
with $c$ a constant independent of $D$.  The result now follows from
(7.2).
\enddemo\sskip

\demo{Proof of Proposition 7.1} We shall first prove that
$\pmb\delta$
induces an isomorphism in cohomology.
Let $\widetilde\bfS_P=(\text{\bf Div}^{r_{l-1},d_{l-1}}_{C/k})^{P_{l-1}}$.
By induction on $l$ it
suffices to prove the result for the map
$$
\pmb\delta: \widetilde\bfS_P \times
(\text{\bf Div}^{r'_l,d'_l}_{C/k})^{ss}
\longto\bfS_P
\tag 7.3
$$
defined by $\pmb\delta ((F,G))=G\oplus F$.
Since the image of $\pmb\delta $ is contained in $\bfS_P^o$,
we can substitute $\bfS_P^o$ for $\bfS_P$ in (7.1).

We shall show that
for every $i\geq 0$, there is an integer $N$ such that given any $D$, with
$\deg D > N$,  one can find open sets $A$ and $V$, of the
elements of index $D$ in the
left and right hand sides of (7.3), so that in the commutative diagram
$$
\CD
H^i(S_P(D))  @>{\pmb\delta^*}>>
H^i (\widetilde S_P(D)\times
(\text{Div}^{r'_l,d'_l}_{C/k}(D))^{ss}) \\
@VVV  @VVV \\
H^i(V) @>{\pmb\delta^*}>> H^i (A)
\endCD
$$
both the vertical and bottom arrows are isomorphisms.
\enddemo

Let $D$ be an effective divisor on $C$ and let
$\text{\bf x}$ be a nonempty set of
points in $C$ not lying on the support of  $D$.
Fix the direct sum decomposition
$$
\bfo_{C}(D)^r=\bfo_{C} (D)^{r'_l}\oplus\bfo_{C}(D)^{r_{l-1}}
$$
and define the following open sets:
\roster
\item  $A^{(1)}_{\x}(D)\subset \widetilde S_P(D)$
consists of the divisors $F$ such that
the support of $\bfo_{C}(D)^{r_{l-1}}/F$ is disjoint from
$\text{\bf x}$.

\item
$A^{(2)}_{\x}(D)\subset
 (\text{Div}^{r'_l,d'_l}_{C/k}(D))^{ss}$
consists of the divisors $G$ such that
the support of $\bfo_{C}(D)^{r'_l}/G$ does not intersect
$\text{\bf x}$.

\item
$V_{\x}(D)\subset
S_P^o(D)$ consists of the divisors $E$,  having
Harder-Narasimhan filtration
$$
0=E_0\subset\dots \subset E_{l-1}\subset E_l=E \ ,
\tag 7.4
$$
such that
the points of $\text{\bf x}$
do not lie on the supports of the quotients
$\bfo_{C}(D)^r/E$ and   $\bfo_{C}(D)^{r_{l-1}}/\ov E
_{l-1}$.  Here $\ov E_{l-1}$  is  the image of
$E_{l-1}$ under the projection of
$\bfo_{C}(D)^r$
into $\bfo_{C}(D)^{r_{l-1}}$.

\item
 $W_{\x}(D)\subset
  (\D^{r_{l-1}, d_{l-1}}_r(D))^{P_{l-1}}$
consists of the submodules $L$ of
$\bfo_{C}(D)^r$ such that
$L$ projects injectively
onto a submodule $\ov L$ of $\bfo_{C}(D)^{r_{l-1}}$ and
the support of
$\bfo_{C}(D)^{r_{l-1}}/\ov L$ does not intersect $\text{\bf x}$.
\endroster
\sskip\sskip

\proclaim{Proposition 7.3}
Let
$$
%N= (2g + k - 1) +\max_{1\leq i\leq\l} {d'_i\over r'_i }\ .
N_q=\max \ \{  \  {d'_1\over r'_1 }+2g -1   , \   {d'_l\over r'_l }+2g -1 +
q \}
\tag 7.5
$$
for $q\ge 0$. Let $\x$ be a nonempty set of points in $C$ and set $N=N_{\#\x}$,
where $\#\x$ is the cardinality of $\x$.
If $D$ is an effective divisor  with $\deg D >N$, then
the open sets just defined have the following properties:
\roster
\item If we set $A_{\x}(D)=A^{(1)}_{\x}(D)\times A^{(2)}_{\x}(D)$, then
$$
\pmb\delta : A_{\x}(D)\longto  V_{\x}(D)
$$
induces an isomorphism of cohomology groups.

\item If we set
$$
A^{\x }_D =\bigcup_{x\in\text{\bf x}}
 A_{\{ x \} }(D)
$$
and
$$
V^{\x}_D = \bigcup_{x\in\text{\bf x}} V_{\{ x \} }(D)\ ,
$$
then
$$
\pmb\delta : A^{\x}_D \longto  V^{\x}_D
$$
also induces  an isomorphism of cohomology groups.

\item  The codimensions of the complements of both
$A^{\x}_D$ and $V^{\x}_D$ are
not less than the integral part of half the cardinality of $\x$.
In particular, for $i$ smaller than this cardinality, the
$i$-th cohomology groups
of both $A^{\x}_D$ and $V^{\x}_D$ are isomorphic, under restriction,
to those of their ambient varieties.
\endroster
\sskip\sskip
It follows that the open sets $A$ and $V$, mentioned  before,
may be obtained by choosing a set $\x$ of more than
$i$  distinct points in $C$, not appearing on $D$, and
setting $A=A^{\x}_D$,  $V=V^{\x}_D$.
\endproclaim
\demo{Proof}
$(1)$\ \  Let $M=M(r'_l,d'_l\,;\,\x )$
be  the moduli space  of $\x$-rigidified semistable bundles,  having
rank $r'_l$ and degree $d'_l$,   that was considered in
the previous section.

Our strategy will be to construct a
commutative diagram
$$
\CD
A^{(1)}_{\x}(D)\times A^{(2)}_{\x}(D) @>\delta>> V_{\x}(D)\\
@V 1\times\varphi VV @VV(\pi,\psi) V\\
A^{(1)}_{\x}(D)\times M
@>> j\times 1> W_{\x}(D)\times M
\endCD
$$
such that the top varieties are affine bundles over the bottom ones
and  $j$ induces an isomorphism in cohomology. ( See, for example,
Proposition 5.12 and Theorem 5.15 in Chapter VI of [23] for the acyclicity of
affine bundles. )

Let us now describe the maps that appear in the diagram above.

Any $G$ in $A^{(2)}_{\x}(D)$ is naturally $\x$-rigidified as follows.
Consider the short exact sequences
$$
0\to\bfo_{C}^{r'_l}\to\bfo_{C}(D)^{r'_l}\to\bfo%{C/k}
^{r'_l}_D\to 0
\tag 7.6
$$
and
$$
0\to G\to\bfo_{C}(D)^{r'_l}\to \bfo_{C}(D)^{r'_l} / G \to  0 \ .
$$
It is clear that their restriction to $\x$ provide us with an
$\x$-rigidification  $u$ determined by the isomorphisms
$$
\bfo_{\x}^{r'_l}\to\bfo_{C}(D)^{r'_l}|_{\x}\leftarrow G|_{\x} \ .
$$
The map $\varphi : A^{(2)}_{\x}(D)\to M$ sends $G$
to $[G,u]$, the isomorphism class  of $(G,u)$.

Similarly, the map $\psi : V_{\x}(D)\to M$ associates to any $E$ in
$V_{\x}(D)$ the isomorphism class of $(E/E_{l-1},u)$ where
$E_{l-1}$ is as in (7.2) and $u$ is the $\x$-rigidification
obtained as follows. The short exact sequences
$$
0\to\bfo_{C}(D)^{r'_l}\to\bfo_{C}(D)^r/
E_{l-1}\to
\bfo_{C}(D)^{r_{l-1}}/\ov E_{l-1}\to 0
$$
and
$$
0\to E/E_{l-1}\to \bfo_{C}(D)^r/{E_{l-1}}\to\bfo_{C}(D)^r/E\to 0
$$
along with (7.6)  provide us, after restriction to $\x$, with the
chain of isomorphisms
$$
\bfo^{r'_l}_{\x} \to \bfo_{C}(D)^{r'_l}|_{\x}  \to
(\bfo_{C}(D)^r/{E_{l-1}})|_{\x} \leftarrow (E/E_{l-1})|_{\x}
$$
which defines $u$.

The  map \ $\pi:V_{\x}(D)\to W_{\x}(D)$
associates to   a divisor $E$   the term
$E_{l-1}$ of its Harder-Narasimhan filtration (7.4).

Finally,
the map \ $j:A^{(2)}_{\x}(D)\to W_{\x}(D)$ \
is defined by \ $j(M)=0\oplus
M$ where
$0\oplus M\ \subset\bfo_{C}(D)^{r'_l}\oplus \bfo_{C}(D)^{r_{l-1}}$.

Using the universal bundles over $W_{\x}(D)$ and $M$
it is possible to construct a vector bundle $\H_{\x} (D)$
over $W_{\x}(D)\times M$
such that  its fibre over $(L,[F,u])$
is canonically isomorphic to $\text{Hom }(F,\bfo_{C}(D)^r/L)$. ( Note
that    $\deg D >N$  implies $\text{Ext}^1(F,\bfo_{C}(D)^r/L)=0$. )

One can lift $(\pi,\psi)$ to a map
$$
\lambda:V_{\x}(D)\to\H_{\x}(D)
$$
which  assigns  to $E$ the homomorphism $E/E_{l-1}\to
\bfo_{C}(D)^r/E_{l-1}$ deduced  from the inclusion
$E\subseteq \bfo_{C}(D)^r$.
The map $\lambda$ is clearly injective.

We claim that $\lambda$ defines
an isomorphism of $V_{\x}(D)$ with
an affine subbundle of $\H_{\x}(D)$.

Let $(L,[F,u])$ be a point of $W_{\x}(D)\times M$
and let  $\mu :F|_{\x}\to G|_{\x}$, $G=\bfo_{C}(D)^r/L$,
be the composition of the chain of isomorphisms
$$
F|_{\x}\buildrel {u^{-1}}\over\to\bfo_{\x}^{r'_l}\to
\bfo_{C}(D)^{r'_l}|_{\x}\to G|_{\x}
$$
Let $i$ be the inclusion of $\x$ in $C$ and let
$$
\mu_0: F\longto i_*i^* G\simeq G\otimes i_*\bfo_{\x}
$$
be the composition of $F\to i_*i^* F$ and $i_*\mu$.

It is easily verified that the elements of
$\text{Hom }(F, G)$ in
the image of $\lambda$ are precisely those which map
to $\mu_0$  under the map
$$
\text{Hom }(F, G)\longto
\text{Hom }(F,  G\otimes i_*\bfo_{\x}) \ .
$$

Since $\deg D > N$ we have
$$
\text{Ext}^1(F, G\otimes\bfo_{C}(-\x )\, )=0
$$
and the sequence
$$
0\to\text{Hom }(F,  G\otimes\bfo_{C}(-\x )) \to\text{Hom }
(F, G)\to\text{Hom }(F, G\otimes i_* \bfo_{\x})\to 0
$$
is exact.

Thus  the image of $\lambda$ intersects every fibre of the vector
bundle $\H_{\x}(D)$ in an affine space of constant dimension
and is therefore an affine bundle over  $W_{\x}(D)\times M$.
It is not difficult to verify
that the differential of $\lambda$ is everywhere injective
( note that the tangent space to
$M$ at $[F,u ]$ can be
identified with $\text{Ext}^1(F,F(-\x))$.)
It follows that
$\lambda$ is a closed immersion and the claim is verified.

In a similar way  one can prove that
$\varphi: A^{(2)}_{\x}(D)\longto M$ is also an  affine bundle.

We now turn  to the proof that $j:A^{(1)}_{\x}(D)\longto W_{\x}(D)$
induces an isomorphism in cohomology.
Consider the retraction $r:W_{\x}(D)\longto A^{(1)}_{\x}(D)$
defined by $r(L)=\ov L$, where $\ov L$ is the projection of
$L\subseteq\bfo_{C}(D)^{r}$ into
 $\bfo_{C}(D)^{r_{l-1}}$.
The fibre of $r$ over $\ov L$ is isomorphic to $\text{Hom }(\ov L,\,
\bfo_{C}(D)^{r'_l})$. Since $\deg D>N$ we have $\text{Ext}^1(\ov L,\,
\bfo_{C}(D)^{r'_l})=0$,  and  $W_{\x}(D)$ is a vector bundle over
$A^{(1)}_{\x}(D)$  having $j$ as its zero section. This proves
what we wanted.

$(2)$\ \  Observe   that if $\y$  is another
nonempty set of distinct points of $C$,
then
$$
A_{\x}(D)\cap A_{\y}(D)=A_{\x\cup\y}(D)
$$
and
$$
V_{\x}(D)\cap V_{\y}(D)=V_{\x\cup\y}(D)\ .
$$
Comparing obvious Mayer-Vietoris sequences and using the
five-lemma, we  conclude that $\pmb\delta : A^{\x}_D \longto  V^{\x}_D$
gives an isomorphism in cohomology.

$(3)$\ \   We shall  only consider the complement of $V^{\x}_D$, but
the same argument applies  to that  of $A^{\x}_D$.
Remark that if $E\in S_P^o(D)- V^{\x}_D $,
then at least half the points of $\x$ belong simultaneously either to
the support of $\bfo_{C}(D)^r/E$ or  to that of
$\bfo_{C}(D)^{r_{l-1}}/\ov E_{l-1}$.
We are thus reduced to proving the following
assertion.

Let $\y$ be a nonempty set of points in $C$. If
$Z_{\y}(D)$ is the subset
of $\div(D)$ consisting of the divisors $E$ such that $\y$ is contained
in the support of
$\bfo_{C}(D)^r/E$, then the codimension of
$Z_{\y}(D)$ is greater, or equal, than  the number of elements in $\y$.

This   can be established by induction on $r$ using the
stratification of $\div(D)$ associated to the direct sum decomposition
$ \bfo_{C}(D)^r= \bfo_{C}(D)\oplus\bfo_{C}(D)^{r-1} $.
\enddemo
\sskip\sskip

To complete the proof of Proposition 7.1 we still have   to show
that the cohomology of $\bfS_P$ stabilizes.

Since $\pmb\delta$ induces an isomorphism
in cohomology, it suffices to show that  the cohomology of
$(\text{\bf Div}^{r,n }_{C/k})^{ss}$
stabilizes.

Given $i\geq 0$, let $N=N_{i+1}$ be as in (7.5).
If $D,D'$ are effective divisors, with $D'\geq D$ and $\deg
D > N$, then choose $\x$ to consist of
$i+1$ distinct points
in $C$ not lying on $D'$. Now, in the
commutative diagram
$$
\matrix
(\text{\bf Div}^{r,n }_{C/k}(D))^{ss}& @>>>  &
(\text{\bf Div}^{r,n }_{C/k}(D'))^{ss} \\
 \uparrow  & &\uparrow  \\
A^{(2)}_{\x}(D)  & @>>>  &  A^{(2)}_{\x}(D')  \\
\varphi\searrow &  &\swarrow\varphi' \\
& {  M } &
\endmatrix
$$
both $\varphi$ and $\varphi '$ are affine bundles over $M$. It follows that
the top arrow
induces an isomorphism of $i$-th cohomology groups, and we are done.

\head 8. The Abel-Jacobi map \endhead

Recall that there is a coarse moduli space
$N(r,n)$ parametrizing isomorphism classes of stable vector
bundles of rank $r$ and degree $n$ over the curve $C$. This is
a smooth quasi-projective variety of dimension $1+r^2(g-1)$. When
$r$ and $n$ are coprime the notion of stable and semistable
bundle coincide and $N(r,n)$ is a smooth projective algebraic variety.
Moreover, in this case, $N(r,n)$ is a fine moduli space. In particular there
are Poincar\'e bundles $\gotP^{r,n}$ over $C\times N(r,n)$ such
that for every $[E]\in N(r,n)$ the restriction $\gotP^{r,n}_{[E]}$
of $\gotP^{r,n}$ to $C\times \{[E]\}$ is isomorphic to $E$.

Let $r$ and $n$  be coprime.
It is natural to define, by analogy with the
classical case, Abel-Jacobi maps
$$
\pmb\vartheta:(\text{\bf Div}^{r,n}_{C/k})^{ss}\longto N(r,n)\ ,
$$
by assigning to a divisor $E$ its isomorphism class as a
vector bundle.
Here $N(r,n)$ is considered as a constant ind-variety.

The aim of this section is to prove the following

\proclaim{Proposition 8.1} If $r$ and $n$ are coprime,
then there is a ``coherent locally free'' module $\pmb E$, defined over
the constant ind-variety $N(r,n)$, and a morphism
$$
\matrix
({\text{\bf Div}}^{r,n}_{C/k})^{ss}& & @>\text{\bf i}>> & &\Bbb P(\pmb E )\\
_{\displaystyle\pmb\vartheta} &\searrow & & \swarrow & \\
&& N(r,n) &&
\endmatrix
$$
inducing an isomorphism in cohomology. In particular,
$$
H^*((\divb)^{ss} ) = H^*(N(r,n))[x]\ ,
$$
where $x$  is an independent variable of degree 2.  This gives an identity of
Poincar\'e series
$$
P(N(r,n);t)=(1-t^2)P((\divb)^{ss};t)\ .
$$
\endproclaim
\demo{Proof}
For simplicity, we shall work over the cofinal subset of positive divisors such
that
$$
\deg D > 2g + {\displaystyle n\over\displaystyle  r}-1 \ .
$$

We construct,
using
the Poincar\'e bundle,  a vector bundle $\H(D)$
over $N(r,n)$ whose fibre at any $[E]$ is canonically isomorphic
to $\text{Hom }(\gotP_{[E]},\bfo_{C}(D)^r)$. ( Recall
that for any stable $(r,n)$-bundle $F$, one has
$\text{Ext}^1(F,\bfo_{C}(D)^r)=0$.~)

There is a natural morphism $i_D:(\div(D)))^{ss}\longto \Bbb P(\H(D))$
defined by sending any $E$ to the class in
$\Bbb P(\text{Hom }
(\gotP_{[E]},\bfo_{C}(D)^r))$ of the composite of any isomorphism
$\gotP_{[E]}\longto E$ with the inclusion $E\rightarrow\bfo_{C}(D)^r$.
Recall that the automorphisms of a stable bundle are given by
multiplication by non-zero scalars.

( We leave to the conscientious reader the task of
defining $\H(D)$, {\bf i}
using E.G.A. III $\S 7.7$ and
the universal properties of the different objects under consideration. )

The map $i_D$ is injective by definition. Also, it is not difficult
to verify that $i_D$ is an immersion (let us  mention here that the
differential of $i_D$ at $E$ can be naturally identified with the
coboundary map $\text{Hom }(E,\bfo_{C}(D)^r/E)\longto\text{Ext}^1
(E, E))$. Furthermore it is clear that the image of $i_D$ intersects
any fibre $\Bbb P(\text{Hom }(\gotP_{[E]},\bfo_{C}(D)^r))$ in
the open subset
consisting of the classes of injective homomorphisms.
\enddemo

We claim that $\text{\bf i}$ is a quasi isomorphism. This follows from
the following
\proclaim{Lemma 8.2} Let $E,F$ be locally free $\bfo_C$-modules
of rank $r$
such that $\text{\rm Ext}^1(E,F)$ vanishes. For any effective divisor $D$ let
$c_D$ be the codimension in $\text{\rm Hom }(E,F(D))$ of the closed subset
consisting of the homomorphisms which are not injective. Then we have
$$
c_D\geq\deg D\ .
$$
\endproclaim
\demo{Proof}
Since $F\otimes\bfo_D$ is isomorphic to $\bfo^r_D$, we have
an exact sequence
$$
0\longto F\longto F(D)\longto\bfo^r_D\longto0 \ .
$$
 From this we obtain another exact sequence
$$
0\longto\text{Hom }(E,F)\longto\text{Hom }(E,F(D))\buildrel\Phi\over
\longto\text{Hom }(E,\bfo^r_D)\longto0\ .
$$
If $\lambda\in\text{Hom }(E,F(D))$ is not injective then
the rank of the image of the homomorphism
$\Phi(\lambda)$ is smaller than $r$, and $\Phi(\lambda)$ factors
through a surjection
$E\longto\bfo^s_D$ for some $s<r$. Thus, in order to establish the lemma,
it suffices to show that the subset of homomorphisms having this
propriety has codimension not less than $\deg D$ in $\text{Hom }(E,
\bfo^r_D)$.

Note that if $D=\sum^q_{i=1}n_i P_i$, then $\bfo^r_D\cong\oplus
^q_{i=1}\bfo^r_{D_i}$, with  $ D_i={n_iP_i}$, and
$$
\text{Hom }(E,\bfo^r_D)\cong\oplus
^q_{i=1}\text{Hom }(E,\bfo^r_{D_i})\ .
$$
Hence we may assume that $D=m\,P$. Let $\Psi^s$
be the map
$$
\text{Hom}^{\text{Surj}}(E,\bfo^s_D)\times\text{Hom }(\bfo^s_D,
\bfo^r_D)\longto\text{Hom }(E,\bfo^r_D)
$$
given by composition. The group $G$ of automorphisms of $\bfo^s_D$,
where $D=mP$,
has dimension $s^2m$ and
acts freely on the domain of $\Psi^s$ by $g(\alpha,\beta)=(g\circ
\alpha,\beta\circ g^{-1})$.
Since $\Psi^s$ is clearly constant along the orbits, we obtain
$$
\dim \text{\rm Im}\ \Psi^s\leq 2rsm-s^2m\
$$
and
$$
\text{codim Im}\, \Psi^s\geq(r-s)^2m\
$$
which proves the lemma.
\enddemo
\head 9. The fundamental action  on the ind-variety of divisors \endhead

Let  $G$ be the group of automorphisms of the constant $\bfo_{C}$-module
$K^r$. Clearly $G$ is isomorphic to $\text{GL}(r,K)$.
We would like to construct an ind-variety $\text{\bf G}$ having $G$ as its
$k$-rational  points. We proceed as follows.

First, for any positive divisor   $D$,   we set
$$
G(D)=\{g\in G\ \vert\  g(\bfo_{C}^r)\subseteq\bfo_{C}(D)^r\}\ .
$$
Clearly
$$
G=\bigcup_{D\geq0} G(D)\ .
$$
Note that the restriction function
$$
G(D) \longto
\text{Hom}^{\text{Inj}}(\bfo_{C}^r,\bfo_{C}(D)^r)
$$
is a bijection (
 an automorphism of $K^r$ is determined by its restriction to $\bfo_{C}^r$
and,  conversely, any injective homomorphism $\bfo_{C}^r\longto K^r$
extends to an automorphism of $K^r$. )
It follows that we can identify $G(D)$ with the $k$-points of an open
subvariety of the affine space $\Bbb V (E_D)$ determined by the $k$-vector
space $E_D=\text{Hom }(\bfo_{C}^r,
\bfo_{C}(D)^r)$ .

If  $D\leq D'$, then the  inclusion $G(D)\hookrightarrow G(D')$ is induced by
a closed immersion of the corresponding varieties. The resulting  inductive
system is $\bfG$.

The  composition law $G\times G\to G$ is induced from a morphism of algebraic
ind-varieties
$$
\pmb\gamma :\bfG\times\bfG\longto\bfG
$$
determined by the obvious morphisms of
algebraic varieties
$$
G( D)\times G( D')\longto G( D+D')\ ,
$$
one for each pair $(D,D')$ of positive divisors.

This makes $\bfG$ into an  ind-variety  in monoids
( note, however,  that if
$i:G\to G$ is the map giving the inverse, then
for any $D$ there is no $D'$
such that $i(G( D))\subseteq G( D')\,$. )

We also define an ind-variety in monoids $\ov{\bfG}$ corresponding to
$\ov G=G/k^*$.
One has
$$
\ov G=\bigcup_{D\geq0}\ov G( D)
$$
where $\ov G( D)=
G(D)/k^*$. Since $\ov G( D)$  can be identified with the
$k$-points of an open subvariety of $\Bbb P (E_D)$, namely the
image of $E^{\text{Inj}}_D=\text{Hom}^{\text{Inj}}(\bfo_{C}^r,
\bfo_{C}(D)^r)$  under  the projection
$$
\Bbb V (E_D)-\{ 0\}\longto \Bbb P (E_D)\ ,
$$
we can take $\ov{\bfG}$ to be the  system determined by these
subvarieties.
Note that if $\text{\bf E}=\{ E_D \}_{D\geq 0}$ is the obvious
module over the constant ind-variety $\text{Spec }k$, then
$\bfG$ ( resp. $\ov{\bfG}$ ) is an open ind-subvariety of
$\bold V=\Bbb V(\text{\bf E})$
( resp. $\bold P=\Bbb P (\text{\bf E} )$.~)

\proclaim{Proposition 9.1}
Let $\pmb\tau$ be the restriction of $\bfo_{\text{\bf P}} (1)$ to the
open ind-subvariety $\ov\bfG$.
The cohomologies  of $\bfG$ and   $\ov{\bfG}$
stabilize. Moreover, if  $t$  is the   Chern class of $\pmb\tau$, then
$$
H^*(\bfG)=\Bbb Z_{\ell}
$$
and
$$
H^*(\ov{\bfG})=\Bbb Z_{\ell}\, [\,t\,]\ ,
$$
with $t$ algebraically  independent over $\Bbb Z_{\ell}$.
\endproclaim

\demo{Proof}  Both the inclusions $\bfG\hookrightarrow\text{\bf V}$ and
$\ov\bfG\hookrightarrow\text{\bf P}$ are quasi-isomorphisms
( proof as in Lemma 8.2. )
\enddemo

There is  a natural   action of $G$ on the set of $(r,n)$-divisors
$\div$ given by $g\cdot E =g(E)$.
The  orbits of  $ \div$ under this action are in one-to-one correspondence
with the set of isomorphism
classes of $(r,n)$-vector bundles   over $C$.  Note that
any such  bundle can be injected into $K^r$ and that any isomorphism between
two $(r,n)$-divisors can be extended to an automorphism of $K^r$.

The action $G\times\div\longto\div$ is induced from an action
$$
\pmb\alpha:\bfG\times\divb\longto\divb\
$$
determined by the obvious
morphisms of algebraic varieties
$$
\alpha_{D,D'}:G(D)\times\div(D')\longto\div (D+D')\ ,
$$
one for each pair $(D,D')$ of positive divisors.
Since $k^*$ acts trivially on $\divb$, we also have
an action
$$
\ov{\pmb\alpha}:\ov{\bfG}\times\divb\longto\divb\ .
$$
\bigskip

Let  $\pmb{\U}\,^{r,n}$  be the universal  $(r,n)$-divisor  over
$C\times\divb$.
This is obtained, as usual, from the universal quotients over the components
of $\divb$.
We shall also denote  $\pmb{\U}\,^{r,n}$ the associated bundle $\Bbb V
(\pmb{\U}\,^{r,n} )$.
This bundle has an obvious $\bfG$-linearization
$$
\tilde{\pmb\alpha}:\bfG\times\pmb{\U}\,^{r,n}\longto\pmb{\U}\,^{r,n}
$$
which is compatible with $\pmb\alpha$.

Let $\bfS$ be an ind-variety having an action
$$
\pmb\beta:\ov{\bfG}\times\bfS\longto\bfS
$$
of $\ov{\bfG}$ on it ( e.g. $\bfS=\bfS_P$ with the action induced by
$\ov{\pmb\alpha}$. )

If $\bold f:\bfS\longto\divb$ is a $\ov{\bfG}$-equivariant
morphism, then we set
$$
\pmb{\V}=(1_C\times \bold f)^*\pmb{\U}\,^{r,n}\ .
$$
Given points
$x_0$ in $C$ and $s_0$ in $\bfS$ respectively, we consider the sections
$$
C\buildrel\pmb\eta_{s_0}\over\longto C\times\bfS
\buildrel\pmb\epsilon_{x_0}\over\longleftarrow\bfS
$$
defined by $\pmb\eta_{s_0}(x)=(x,s_0)$ and $\pmb\epsilon_{x_0}(s)=(x_0,s)$.

Also let
$$
\ov{\bfG}\buildrel\bold q_1\over\longleftarrow C\times\ov{\bfG}\times\bfS
\buildrel\bold q_2\over\longto C\times\bfS \ \ ,
$$

$$
\ov{\bfG}\buildrel\bold p_1\over\longleftarrow\ov{\bfG}\times \bfS
\buildrel\bold p_2\over\longto\bfS
$$
and
$$
C\buildrel\bold u_1\over\longleftarrow C\times\ov{\bfG}
\buildrel\bold u_2\over\longto \ov{\bfG}
$$
be the natural projections.

\proclaim{Lemma 9.2} With the above notations, we have:
\roster
\item $(1_C\times\pmb\beta)^*\,\pmb{\V}\
=\ \bold q^*_1\,\pmb\tau \otimes \,\bold q^*_2\,\pmb{\V}$

\item $\pmb \beta\,^*\;\pmb\epsilon_{x_0}^*\pmb{\V}\
=\ \bold p_1^*\pmb\tau\otimes \,\bold p^*_2\,\pmb\epsilon_{x_0}^*\pmb{\V}$

\item $(1_C\times\, (\pmb\beta\circ\pmb\eta_{s_0}))^*\,\pmb{\V}\
=\ \bold u^*_2\,\pmb\tau \otimes
\,\bold u^*_1\,\pmb\eta_{s_0}^* \pmb{\V}$.
\endroster
\endproclaim
\medskip

\demo{Proof} $(2)$ and $(3)$ are immediate consequences of $(1)$.
By functoriality,
it suffices to  consider the case where $\bold S=\divb$ and
$\bold f $ is the identity morphism. For $D\leq D'$,  there is
a commutative diagram
$$
\CD
G( D)\times\U^{r,n}(D') @>\tilde\alpha_{D,D'}>> \U^{r,n}(D+D')\\
@V VV @V VV \\
C\times\ov G(D)\times\div(D') @>> 1_C\times\ov\alpha_{D,D'}>
C\times \div(D+D')\ .
\endCD
$$
and we can identify $G(D)\times\U^{r,n}(D')$ with the complement of the
zero section in the bundle
$(q_1)_{D,D'}^*\tau(D)\otimes (q_2)_{D,D'}^*\U^{r,n}(D')$.
Moreover $\tilde\alpha_{D,D'}$ extends across the zero section to give
an isomorphism
$$
(q_1)_{D,D'}^*\tau (D)\otimes (q_2)_{D,D'}^*\U^{r,n}(D')
\longto (1_C\times\ov\alpha_{D,D'})^*\U^{r,n}(D+D')
$$
and  the lemma follows.
\enddemo

\proclaim{Proposition 9.3}
In the situation of the previous lemma let
$\ell$ be prime with $r$ ( or, if $k=\Bbb C$ and
singular cohomology is being used, assume that the coefficient
field is
of characteristic prime with $r$.)
Then there exists an element $x$ in $H^2(\bfS)$ such that
$$
\pmb\beta\,^*(x)=\bold p_1^*(t)+\bold p_2^*(x)\ .
$$
\endproclaim

\demo{Proof}
As in the previous proof we can assume that $\bfS=\divb$. Let $y$ be
the first Chern class of $\pmb\epsilon^*_{x_0}\pmb{\U}\,^{r,n}$.
Then using 2), in the previous lemma,
we have
$$
\ov{\pmb\alpha}\,^*(y)= r\cdot\bold p_1^*(t)+\bold p_2^*(y)\ .
$$
Now take $x={\displaystyle 1\over\displaystyle r}\cdot y$,
and we are done.
\enddemo

\definition{Definition 9.4}
Let $\pmb\beta:\ov{\bfG}\times\bfS\longto\bfS$ be an action
of $\ov{\bfG}$ on the ind-variety $\bfS$.
The  ring of invariants in $H^*(\bfS)$ is defined  by
$$
H^*(\bfS)^{\ov{\bfG}}
=\{\xi\in H^*(\bfS)\,\,\vert\,\,\pmb\beta\,^*(\xi)=\bold p_2^*(\xi )\}\ .
$$
\enddefinition

\proclaim{Proposition 9.5}
In the situation of the previous definition let
$x\in H^2(\bfS)$ be such that
$$
\pmb\beta\,^*(x)=\bold p_1^*(t)+\bold p_2^*(x)\ .
$$
Then we have
$$
H^*(\bfS)=H^*(\bfS)^{\ov{\bfG}}[\,x\,]
\tag{9.1}
$$
and $x$ is
algebraically  free  over $H^*(\bfS)^{\ov{\bfG}}$.
\endproclaim

\demo{Proof}
Identify
$H^*(\ov{\bfG}\times \bfS)$
with
$H^*(\bfS)[\,t\,]$
( note that we  are using the same letter
for both $t$ and $\bold p_1^*(t)$ ),
and let $\mu:H^*(\bfS)[\,t\,]\longto H^*(\bfS)$ be the
homogeneous homomorphism which is the identity on $H^*(\bfS )$ and
which sends $t$ to $-x$.
Define $\varphi:H^*(\bfS)\longto H^*(\bfS)$ as the composite
$\varphi=\mu\circ\pmb\beta\,^*$.
Explicitly if
$\pmb\beta\,^*(\xi)=\sum_{i\geq0}\xi_it^i$,
with $\xi_i\in H^*(\bfS)$, then
$$
\varphi(\xi)=\sum\xi_i(-x)^i\ .
$$
We claim that the image of $\varphi$ coincides with the ring of invariants.
Since $\varphi$ is clearly  the identity on
$H^*(\bfS)^{\ov{\bfG}}$, it suffices to show that the image under $\varphi$
of any $\xi$ is invariant. We shall see below   that  this  follows
from the identity,
in  $H^*(\bfS)[t,u]$,
$$
\sum_{i\geq0}\xi_i(t+u)^i=\sum_{i\geq0}\;\pmb\beta\,^*(\xi_i)u^i\ .
$$
which, in this context, is a consequence of
the associativity of the action  together with the identity
$\ov{\pmb\gamma}^*(t)=t+u$ where $\ov{\pmb\gamma}$ is the composition law on
$\ov\bfG$.
( Note that $\ov{\pmb\gamma}^*(t)$ has to be a linear combination
of $t$ and $u$; an easy calculation shows that  both coefficients are 1. )

Now, as promised,  we have
$$
\eqalign{
\pmb\beta\,^*(\varphi(\xi))
&=\sum_{i\geq0}\pmb\beta\,^*(\xi_i)\pmb\beta\,^*(-x)^i           =\cr
&=\sum_{i\geq0}\pmb\beta\,^*(\xi_i)(-x-t)^i                      =\cr
&=\sum_{i\geq0}\xi_i(-x)^i                                       =\cr
&=\bold p_2^*(\varphi(\xi))\ .
\cr}
$$

In order to prove (9.1) it suffices to show
that every homogeneous element $\xi$ of $H^*(\bfS)$ belongs to
$R=H(\bfS)^{\ov{\bfG}}[x]$. We proceed by induction on $\deg\xi$.
We know that
$\varphi(\xi)=\sum_{i\geq0}\xi_i(-x)^i$ is invariant and,
by induction hypothesis, that every
$\xi_i $,  $i>0$, is already in $R$.
It follows that
 $\xi_0 $ also belongs to $R$. But $\xi_0=\xi$, and  hence $H^*(\bfS )$
coincides with $R$ ( note that $\xi_0=\xi$ because
$s\mapsto (1_{\ov{\bfG}},s)$ defines a section of $\pmb\beta$. )
Finally,
since
$\pmb\beta\,^*(x)=x+t$ with $t$ algebraically independent over $H^*(\bfS)$,
the element $x$ is algebraically independent
over the ring of invariants,
and we are done.
\enddemo
\sskip
\proclaim{Corollary 9.6}  In the situation of  Proposition 9.3 one has
$$
H^*(\bfS)=H^*(\bfS)^{\ov{\bfG}} [\,x\,]\ .
$$
\endproclaim

\medskip
\remark{Remarks}

$(1)$\ \  In the situation of  Proposition 9.5 the morphism
$(\pmb\beta\circ\pmb\eta_{s_0})^*$,   $s_0\in\bfS$,
sends $x$ to $t$ and coincides on $H^*(\bfS)^{\ov{\bfG}}$
with the map induced in cohomology by the inclusion $\{s_0\}\hookrightarrow
\bfS$.

$(2)$\ \  The corollary applies in particular to the Shatz
strata $\bfS_P$  and gives  a generalization of the isomorphism
$$
H^*((\text{\bf Div}^{r,n}_{C/k})^{ss})\simeq H^*(N(r,n))[\,x\,]
$$
that was  established in Section 7  under the assumption that
$r$ and $n$  be coprime.
In this case, the ring of invariants coincides with the image of
$H^*(N(r,n)$ under $\pmb\vartheta^*$.

$(3)$\ \ The assumption made in Proposition 9.3 about the prime
$\ell$  can be relaxed but not eliminated;
it is not difficult to describe an explicit set of generators
of $H^2(\divb)$, and to verify that the conclusion of the lemma holds
if,  and only if, $\ell$ is prime  to the greatest common denominator
of $n$ and $r$.
\endremark

\head 10. Perfection of the Shatz stratification \endhead

Recall that for every effective divisor $D$ we have
 the Shatz stratification\
$$
\div(D)=\bigcup_{P\in\Cal P_{r,n}}
S_P(D)\ ,
$$
where $S_P(D)$ is
a locally closed subscheme,
consisting of the divisors having Harder-Narasimhan filtration of type $P$.

If $\deg D$ is large enough, then $S_P(D)$ is smooth and has codimension
$\d_P$ independent of $D$. As $D$ varies the varieties $S_P(D)$
form an inductive system  $\bfS_P$ and one has
the stratification
$$
\divb=\bigcup_{P\in\Cal P_{r,n}}\bfS_P\ .
$$

The set $\Cal P=\Cal P_{r,n}$  of Shatz $(r,n)$-polygons
( i.e.  strictly convex polygons in $\Bbb R^2$
joining $(0,0)$ with $(r,n)$ ) is partially ordered by the relation
$P\leq P'$ if, and only if,  $P'$ lies above $P$.
A subset $I$ of $\Cal P $ is open if  $P\in I$
and $P\leq P'$ imply $P'\in I$. If $I$ is open, then
$$
S_I(D)=\bigcup_{P\in I}
S_P(D)
$$
is an open subset of $\div (D)$ for each positive $D$
( same statement for $\bold S_I=\bigcup_{P\in I}\bfS_P$. )

Let $I\subset\Cal P$ be   open  and let  $P$ be a minimal element  in the
complement of $I$. The set $J=I\cup \{P\}$ is still open and, as a
consequence of Proposition 4.1,
$S_P(D)$ is closed in the open set $S_J(D)$.

Now if $D\leq D'$, with $\deg D $ large enough, then we  have
a commutative diagram
$$
\matrix
&\to& H^{i-2\d_P}(S_P(D'))&\longto &H^i(S_J(D'))&
\longto &H^i(S_I(D'))&
\to&\cr
&&\downarrow&&\downarrow&&\downarrow&&\cr
&\to& H^{i-2\d_P}(S_P(D))&\longto &H^i(S_J(D))&\longto &H^i(S_I(D))&
\to&
\endmatrix
$$
where the horizontal arrows are Gysin sequences.
Using this diagram and the five-lemma one can show by induction on the
cardinality of $I$ that  the cohomology of $H^*(\bfS_I)$, $I$ finite,
stabilizes. Moreover, taking the projective limit
of these  sequences  we obtain the Gysin sequence
$$
\to H^{i-2\d_P}(\bfS_P)\longto H^i(\bfS_J)\longto H^i(\bfS_I)\to
\tag 10.1
$$
which is still exact.

Throughout the rest of this section, $\ell$-adic  cohomology will
be $\Bbb Q_{\ell}$-cohomology.
We shall say that the stratification
is $\Bbb Q_{\ell}$-perfect, or perfect for short,  if for every finite
open subset $I$ of $\Cal P $ and for
every one of the finitely many minimal elements $\lambda$
of its complement, the  Gysin sequence
(10.1) splits into short exact sequences
$$
0\to H^{i-2\d_P}(\bfS_P)\longto H^i(\bfS_J)\longto H^i(\bfS_I)\to
0\ .
$$

Recall from Proposition 5.1 that for every $d\geq 0$ the number of
strata having codimension smaller than $d$ is finite and that this number
is independent of $D$, for $\deg D$ large enough.

\proclaim{Proposition 10.1}
The Shatz stratification of the ind-variety of $(r,n)$-divisors
is perfect. In particular, there is  an identity of Poincar\'e series
$$
P(\divb;t)=\sum_{P\in\Cal P_{r,n}} P(\bfS_P;t)\ t^{2\d_P}\ .
$$
\endproclaim

\demo{Proof} We use ideas inspired by [2].
As is well known, the composition of the Gysin map
\ $H^{i-2\d_P}(S_P(D))\longto H^i(S_J(D))$
\ with the restriction map
\  $H^i(S_J(D))\longto H^i(S_P (D))$
\ is multiplication by the top Chern class of the normal bundle
$N_P(D)$ of $S_P(D)$ in $\div(D)$.
If $D\leq D' $  and $\deg D$ is large enough,
then the restriction of $N_P(D')$ to $S_P(D)$ is $N_P(D)$.
Therefore the bundles $N_P(D)$, as $D$ varies, define a vector bundle
$\bfN_P$ over $\bfS_P$ ( in fact, over
  a cofinal
subsystem of $\bfS_P$. )
If the top Chern class $e(\bfN_P)$ of $\bfN_P$ is a non zero
divisor in $H^*(\bfS_P)$, then the  Gysin sequence (10.2) will split into
short exact sequences, as desired.

It suffices to prove, with the notation of Proposition 7.1,
that $e(\pmb\delta^*\bfN_p)=\pmb\delta^*(e(\bfN_p))$
is not a zero divisor in
$H^*((\text{\bf Div}^{r'_1,d'_1}_{C/k})^{ss}
\times\dots\times (\text{\bf Div}^{r'_l,d'_l}_{C/k})^{ss})$.

Let $\bfG_i$ and $\ov{\bfG}_i$ be the ind-varieties in
monoids associated to $GL(r'_i,K)$ as
in Section 9.
For $i=1,\dots, l$ fix divisors
$F_i\in(\text{\bf Div}^{r'_i,d'_i}_{C/k})^{ss}$
and consider the morphism
$$
\bold j:\ov{\bfG}_1\times\dots\times \ov{\bfG}_l\longto(\text{\bf Div}^{r'_1,
d'_1}_{C/k})^{ss}\times\dots\times(\text{\bf Div}^{r'_l,d'_l}_{C/k})^{ss}
\tag 10.2
$$
given by $\bold j\;(( g_1,\dots,  g_l))
=( g_1(F_1),\dots, g_l(F_l))$.

Using the K\"unneth formula and the results of Section 8 we can identify
the cohomology ring of the right hand side of (10.2) with
a polynomial ring $R^*[x_1,\dots, x_l]$
where
$$
R^*=\mathop{\otimes}\limits_{1\leq i\leq l}
H^*((\text{\bf Div}^{r'_i,d'_i}_{C/k})^{ss})
^{\ov{\bfG}_i}
$$
Also we have
$$
H^*(\ov{\bfG}_1\times\dots\times \ov{\bfG}_l)
=\Bbb Q_{\ell}[t_1,\dots, t_r]
$$
where
$t_1,\dots, t_r$ are  algebraically independent.
In fact each $t_i$
is the  Chern class of the line bundle $\pmb\tau_i$ obtained
by pulling back to $\ov{\bfG}_1\times\dots\times\ov{\bfG}_l$
the  fundamental line bundle $\pmb\tau$ over $\ov{\bfG}_i$.
Note that the restriction
of $\bold j^*$ to $R^*=\oplus_{i\geq0}R^i$ is
just projection onto $R^0=\Bbb Q_{\ell}$,
and that $\bold j^*$ assigns $t_i$ to each $x_i$.

Using this description of $\bold j^*$ it is immediately verified
that a cohomology class $\alpha$ cannot be
a zero divisor whenever  $\bold j^*(\alpha)\not=0$  ( cf.
$\S 13$ in [2]. )
Hence, in order to
reach our goal we must show that
$$
\bold j^*(e(\pmb\delta^*\bfN_P))=e(\bold j^*\pmb\delta^*\bfN_P)\not=0
$$
Recall from Section 4 that
$$
\bfN_P=R^1 \bold q_*\Cal H om_+(\pmb{\U}_P,\pmb{\U}_P)
$$
where \ $\pmb{\U}_P$ is the universal family
of divisors parametrized by $\bfS_P$, $\Cal H om_+ $
refers to the universal Harder-Narasimhan filtration
of $\pmb{\U}_P$ and $\bold q:C\times \bfS_P\longto\bfS_P$
is the natural projection.

Recall that if   $X$ is a scheme and $\F$ is an $\bfo_X$-module  having a
filtration
$$
0\subseteq\F_1\subseteq\F_2\dots\subseteq\F_l=\F
$$
such that  $\F/\F_1$ is locally free, then  there is a short
exact sequence
$$
0\to\Cal H om_+(\F/\F_1,\,\F/\F_1)\longto\Cal H om_+
(\F,\F)\longto\Cal H om(\F_1,\F/\F_1)\to 0\ .
$$
Using this fact,  it is easily verified that $\Cal H om_+(\pmb{\U}_P,
\pmb{\U}_P)$ is flat over $\bfS_P$.
We also have
$H^2(\Cal H om_+(\U_p(D),\ \U_p(D))\vert_{C\times \{E\}})=0$
for every $E\in S_P(D)$,  and therefore base change commutes
with taking $R^1\bold q_*$.
Hence
$$
\pmb\delta^*(\bfN_P)=R^1\bold q'_*((1_C\times \pmb\delta)^*\Cal H om_+
(\pmb{\U}_P,\pmb{\U}_P))
$$
where
$$
\bold q':C\times \prod_{1\leq i\leq l}(\text{\bf Div}^{r'_i,d'_i}_{C/k})^{ss}
\longto \prod_{1\leq i\leq l}(\text{\bf Div}^{r'_i,d'_i}_{C/k})^{ss}
$$
is the natural projection.

Let $\pmb{\U}_i$ be the pull-back to
$C\times \prod_{1\leq i\leq l}(\text{\bf Div}^{r'_i,d'_i}_{C/k})^{ss}$
of the
universal family of divisors over
$C\times(\text{\bf Div}^{r'_i,d'_i}_{C/k})^{ss}$.
Then we have
$$
(1_C\times\pmb\delta)^*\Cal H om (\pmb{\U}_P,
\pmb{\U}_P)=\mathop{\oplus}\limits_{i,j} \Cal H om (\pmb{\U}_i,\pmb{\U}_j)\ ,
$$

$$
(1_C\times\pmb\delta)^*\Cal H om_-(\pmb{\U}_P,\pmb{\U}_P)=\mathop{\oplus}
\limits_{i\geq j}\Cal H om (\pmb{\U}_i,\pmb{\U}_j)
$$
and
$$
\eqalign{
(1_C\times\pmb\delta)^*\Cal H om_+(\pmb{\U}_P,\pmb{\U}_P) &=
\mathop{\oplus}\limits_{i<j}\Cal H om (\pmb{\U}_i,\pmb{\U}_j) \cr
&=\mathop{\oplus}\limits_{i<j}
\pmb{\U}^\vee_i\otimes\pmb{\U}_j\ .}
$$
It follows that
$$
\pmb\delta^*\bfN_P=\mathop{\oplus}\limits_{i<j}
R^1\bold q'_*\,\pmb{\U}_i^\vee\otimes\pmb{\U}_j\ .
$$
Since  $\bold j^*$ commutes with
$R^1\bold q'_*$, using   Lemma 9.2, we obtain
$$
\bold j^*\pmb\delta^*\bfN_P=\mathop{\oplus}
\limits_{i<j}R^1 \bold p _{2*}(\bold
p^*_1(F_i)^\vee\otimes \bold p^*_2
(\pmb\tau_i)^\vee\otimes\bold  p^*_1(F_j)\otimes\bold  p^*_2(\pmb\tau_j))
$$
where the $F_i$'s are as in (10.2) and
$$
C\buildrel\bold p_1\over\longleftarrow
C\times\ov{\bfG}_1\times\dots\times\ov{\bfG}_l\buildrel\bold p_2\over\longto
\ov{\bfG}_1\times\dots\times\ov{\bfG}_l
$$
are the natural projections. Now, using the projection
formula,  we obtain
$$
\bold j^*\pmb\delta^*\bfN_P=\mathop{\oplus}\limits_{i<j}
(\pmb\tau^\vee_i\otimes\pmb\tau_j)
^{h^1(F^\vee_i\otimes F_j)}\ .
$$
and, therefore,
$$
e(\bold j^* \pmb\delta^*\bfN_P)=\prod\limits_{i<j}(t_j-t_i)
^{h^1(F^\vee_i\otimes F_j)}\not= 0
$$
as required.
\enddemo

\bigskip\bigskip
\np{\bf References}\medskip

\item{[1]}
{\smc E. Arbarello, M. Cornalba, P.A. Griffiths} and {\smc
J. Harris}, {\it Geometry of Algebraic Curves I}, Grundlehren der Math.
Wissenschafter {\bf 267}, Springer-Verlag, New York 1985.
\item{[2]} {\smc
M.F. Atiyah} and {\smc R. Bott}, The Yang-Mills equations over
Riemann surfaces, {\it Phil. Trans. Roy. Soc. London, Ser.} A {\bf 308}
(1982), 523--615.
\item{[3]} {\smc
E. Bifet}, Sur les points fixes du sch\'ema $\text{Quot}_{\Cal O^r_X/X/k}$
sous l'action du tore $\text{\bf G}^r_{m,k}$, {\it C.R. Acad Sci Paris},
t. {\bf 309} (1989), 609--612.
\item{[4]} {\smc
E. Bifet, F. Ghione} and {\smc M. Letizia},  On the Betti numbers of
the moduli space of stable bundles of rank two on a curve,
 Proc. Conf. Vector Bundles and
Special Projective Embeddings, Bergen, July 1989.
\item{[5]} {\smc A. Bialynicki-Birula},  Some theorems on actions of algebraic
groups, {\it Ann. of Math} {\bf 98} (1973), 480--497.
\item{[6]} {\smc
G.D. Daskalopoulos}, {\it The topology of the space of stable bundles
over a compact Riemann surface}, Ph. D. Thesis, University of Chicago.
\item{[7]} {\smc
P. Deligne}, Cohomologie \`a supports propres, Expos\'e XVII, in
{\it Th\'eorie des Topos et Cohomologie Etale des Sch\'emas (SGA $4$) Tome
$3$},
Lecture Notes in Math. {\bf 305}, Springer-Verlag 1973.
\item{[8]} {\smc
P. Deligne}, Th\'eor\`eme de finitude
en cohomologie $\ell$-adique', in
{\it Cohomologie \'Etale (SGA $4\frac 12$)}, Lecture Notes in Math.
{\bf 569}, Springer-Verlag 1977.
\item{[9]} {\smc
U.V. Desale} and {\smc S. Ramanan}, Poincar\'e polynomials of the
variety of stable bundles, {\it Math. Ann.}, {\bf 216} (1975), 233--244.
\item{[10]} {\smc
A. Dold} and {\smc R. Thom}, Une g\'en\'eralisation de la notion
d'espace fibr\'e. Application aux produits symm\'etriques infinis,
{\it C.R. Acad Sci. Paris}, t. {\bf 242} (1956), 1680--1682.
\item{[11]} {\smc
A. Dold} and {\smc R. Thom}, Quasifaserungen und unendliche symmetrische
Produkte, {\it Ann. of Math.}, {\bf 67} (1958), 239--281.
\item{[12]} {\smc
F. Ghione} and {\smc M. Letizia}, Effective divisors of higher rank on
a curve and the Siegel formula, preprint.
\item{[13]} {\smc
A. Grothendieck}, Sur le m\'emoire de Weil [4] G\'en\'eralisation
des fonctions ab\'eliennes,
S\'eminaire Bourbaki 1956/57 no. 141.
\item{[14]} {\smc
A. Grothendieck}, Techniques de construction et th\'eor\`emes d'exis\-tence
en g\'eom\'etrie alg\'ebrique IV: Les sch\'emas de Hilbert,
S\'eminaire Bourbaki 1960/61 no. 221.
\item{[15]} {\smc
G. Harder}, Eine Bemerkung zu einer Arbeit von P.E. Newstead,
{\it J. f\"ur Math}, {\bf 242} (1970), 16--25.
\item{[16]} {\smc
G. Harder}, Semisimple group schemes over curves and automorphic
functions, {\it Actes, Congr\`es intern. Math. 1970, tome 2}, 307--312.
\item{[17]} {\smc
G. Harder} and {\smc M.S. Narasimhan}, On the cohomology groups
of moduli spaces of vector bundles on curves, {\it Math. Ann.}, {\bf 212}
(1975), 215--248.
\item{[18]} {\smc
F. Kirwan}, On spaces of maps from Riemann surfaces to Grassmannians
and applications to the cohomology of moduli of vector bundles,
{\it Ark. f\"or Mat.}, {\bf 24} (1986), 221--275.
\item{[19]} {\smc
G. Laumon}, Correspondence de Langlands g\'eom\'etrique pour les corps
de fonctions, {\it Duke Math. Jour.}, {\bf 54} (1987), 309--359.
\item{[20]} {\smc
H.B. Lawson}, Algebraic cycles and homotopy theory, {\it Ann. of Math.},
{\bf 129} (1989), 319--343.
\item{[21]} {\smc
I.G. Macdonald}, Symmetric products of an algebraic curve, {\it Topology},
{\bf 1} (1962), 319--343.
\item{[22]} {\smc
I.G. Macdonald}, The Poincar\'e polynomial of a symmetric product,
{\it Proc. Camb. Phil. Soc.}, {\bf 58} (1962), 563--568.
\item{[23]} {\smc
J. Milne}, {\it \'Etale Cohomology}, Princeton U. P., 1980.
\item{[24]} {\smc
P.E. Newstead}, Topological properties of some spaces of stable
bundles, {\it Topology}, {\bf 6} (1967), 241--262.
\item{[25]} {\smc
P.E. Newstead}, {\it Introduction to moduli problems and orbit spaces},
Tata Inst. Lecture Notes, Springer-Verlag 1978.
\item{[26]} {\smc
C.S. Seshadri}, {\it Fibr\'es vectoriels sur les courbes alg\'ebriques},
Ast\'erisque 96, 1982.
\item{[27]} {\smc
I.R. Shafarevich}, On some infinite-dimensional groups, {\it Rend.
Mat. Appl. V Ser.}, {\bf 25} (1967), 208--212.
\item{[28]} {\smc
I.R. Shafarevich}, On some infinite-dimensional groups II, {\it Math.
USSR, Izv.}, {\bf 18} (1982), 185--194.
\item{[29]} {\smc
S.S. Shatz}, The decomposition and specialization of algebraic
families of vector bundles, {\it Compositio Math.}, {\bf 35} (1977),
163--187.
\item{[30]} {\smc
J.-L. Verdier} and  {\smc J. Le Potier}, eds. {\it Module des Fibr\'es
Stables sur les Courbes Alg\'ebriques}, Progress in Math. 54, Brikh\"auser,
Boston 1985.
\item{[31]} {\smc
A. Weil}, G\'en\'eralization des fonctions Ab\'eliennes, {\it J. Math.
Pures Appl.}, {\bf 17} (1938), 47--87.

\enddocument
\end